\documentclass[12pt,oneside]{article}
\usepackage{amsmath,amssymb}
\usepackage[all]{xy}
\def\backddots{\mathinner{\mkern1mu\raise1pt
    \vbox{\kern7pt\hbox{.}}\mkern2mu
    \raise4pt\hbox{.}\mkern2mu\raise7pt\hbox{.}\mkern1mu}}
\CompileMatrices
\def\overset#1#2{{\mathop{\kern0pt#2}\limits^{#1}}}
\def\A{{\mathbb A}}

\def\D{{\mathbb D}}
\def\I{{\mathbb I}}
\def\Lm{{\mathbb L}}
\def\R{{\mathbb R}}
\def\Sb{{\mathbb S}}
\def\T{{\mathbb T}}

\def\Adiff{\A-{\rm diff}}
\def\Alin{\A-{\rm lin}}
\def\Ann{\mathop{\rm Ann}\,}
\def\rank{\mathop{\rm rank}\,}
\def\im{\mathop{\rm im}\,}
\def\Ao{\overset{\circ}{\mathbb A}\relax{}}

\def\Fs{\overset{\,\,*}{F}\relax{}}
\def\po{\ea{}^0}
\def\Xo{\overset{\circ}{X}}
\def\Yo{\overset{\circ}{Y}}
\def\ea{\overset{*}{e}{}}

\def\qw{\widetilde{\kern2pt q \kern 2pt}{}}
\newenvironment{definition}
{\par\smallskip {\bf Definition.}~} {\par\smallskip}
\newcounter{examp}[section]
\newenvironment{example}
{\par\medskip \refstepcounter{examp} {\it Example
\thesection.\arabic{examp}.}~} {\par\medskip}

\newcounter{th}[section]
\newenvironment{theorem}
{\par\bigskip \refstepcounter{th} {\bf Theorem
\thesection.\arabic{th}.}~\it} {\par\bigskip}

\newcounter{prop}[section]
\newenvironment{proposition}
{\par\medskip \refstepcounter{prop} {\bf Proposition
\thesection.\arabic{prop}.}~\it} {\par\medskip}

\newcounter{coro}[section]
\newenvironment{corollary}
{\par\medskip \refstepcounter{coro} {\bf Corollary
\thesection.\arabic{coro}.}~\it} {\par\medskip}

\newcounter{lem}[section]
\newcounter{rem}[section]
\newenvironment{remark}
{\par\medskip \refstepcounter{rem} {\bf Remark
\thesection.\arabic{rem}.}} {\par\medskip}

\newenvironment{Proof}
{{\bf Proof.}}{\nopagebreak$\Box$\par\medskip}
\begin{document}
%

\begin{center}
\Large\bf
Lifts of Poisson structures to Weil bundles
\end{center}
\begin{center}
\sc Vadim V. Shurygin, jr
\end{center}

\bigskip

{\sc Abstract.}
  In the present paper, we study complete and vertical lifts
of  tensor fields from a smooth manifold~$M$ to its Weil bundle~$T^\A M$
defined by  a Frobenius Weil algebra~$\A$.
  For a Poisson manifold $(M,w)$, we show that the complete lift $w^C$
and the vertical lift $w^V$ of the Poisson tensor $w$ are Poisson tensors on
$T^\A M$ and establish their properties.
   We prove that  the complete and the vertical lifts induce homomorphisms of the
Poisson cohomology spaces.
  We compute the modular classes of the lifts of Poisson structures.



%
%

\begin{center}
\sc Contents
\end{center}

1. Introduction \hfill \pageref{sect-intro}

2. Poisson manifolds \hfill \pageref{sect-pm}

3. Weil  algebras and Weil bundles \hfill \pageref{sect-wa}

3.1. Weil bundles as smooth manifolds over Weil algebras \hfill
\pageref{subsec-wb}

3.2. The structure of a Frobenius Weil algebra \hfill
\pageref{subsec-struc-fwa}

4. Lifts of tensor fields to Weil bundles \hfill \pageref{sect-lifts}

4.1. Realizations of tensor operations \hfill \pageref{subsec-real}

4.2. The complete lift of a covariant tensor field
\hfill \pageref{subsec-cl-cov}

4.3. The complete lift of a contravariant tensor field \hfill
\pageref{subsec-cl-contrav}

4.4. The vertical lift of a tensor field \hfill \pageref{subsec-vl}

5. Weil bundles of Poisson manifolds \hfill \pageref{sect-wbpm}

5.1. The complete lift of a Poisson tensor \hfill
\pageref{subsec-cl-ps}

5.2. The vertical lift of a Poisson tensor \hfill
\pageref{subsec-vl-ps}

5.3. Modular classes of lifts of Poisson structures \hfill \pageref{subsec-mod}

References \hfill \pageref{sect-ref}

%
%

\section{Introduction}
\label{sect-intro}

  Complete and vertical lifts of Poisson structures from a smooth manifold $M$ to its
tangent bundle $TM$  were studied in several papers, see, e.g.,  J.\,Grabowski
and P.\,Urba\'nski~\cite{G-U, G-U2}, G.\,Mitric and I.\,Vaisman~\cite{M-V,
Vai4}.
  In the present paper we discuss the generalization of these lifts to the case of
a Weil bundle $T^\A M$ for a Frobenius Weil algebra $\A$.
\par
  Differential-geometrical properties of lifts of tensor fields and connections to tangent bundles
were studied by K.\,Yano and S.\,Ishihara~\cite{Yano}.
  Various aspects of geometry of Weil bundles  and, in particular,
lifts of geometric structures were studied by  I.\,Kol\'a\v{r}~\cite{Kolar},
A.\,Morimoto~\cite{Mor}, E.\,Okassa~\cite{Okas2}, L.-N.\,Patterson~\cite{Pat},
P.C.\,Yuen~\cite{Yuen2} and other researchers (see for references the book of
I.\,Kol\'a\v r, P.W.\,Michor and J.\,Slov\'ak~\cite{KMS}).
  Weil bundles were also studied in connection with product preserving
functors~\cite{KMS}.
\par
  A.P.\,Shirokov~\cite{Shir12} discovered that
the Weil bundle $T^\A M$ carries a structure of a smooth manifold over~$\A$,
which made it possible to apply to the study of Weil bundles the theory of
smooth manifolds over algebras  (see, e.g.,~\cite{Shir12, VSh-JMS}).
  In particular,  the theory of realization of
tensor fields and tensor operations over finite-dimensional modules over
Frobenius algebras developed in the papers of
V.V.\,Vish\-nev\-skii~\cite{Vish-74} and G.I.\,Kruch\-kovich~\cite{Kruch}
allowed ones to simplify the theory of lifts of tensor fields and linear
connections from a manifold~$M$ to its Weil bundle $T^\A M$.

This paper is structured as follows.

   In Section 2  we recall basic facts and definitions concerning
Poisson manifolds and Poisson cohomology.

   In Section 3 we recall the notion of a Weil algebra (a local algebra in the sense of
A.\,Weil), basic concepts of the theory of smooth mappings over
Weil algebras
and describe the structure of a smooth manifold over algebra $\A$ on the Weil
bundle $T^\A M$.
  In this section, we also analyze the structure of a Frobenius Weil algebra
and prove some auxiliary statements which we use later.

   In Section 4, we
develop the theory of lifts of skew-symmetric covariant and contra\-variant
tensor fields from a smooth manifold $M$ to its Weil bundle  $T^\A M$ on the
base of the theory of realization of tensor operations.
   We show that the complete lift of exterior forms
induces a homomorphism of the de Rham cohomology spaces $H^*_{dR}(M)\to
H^*_{dR}(T^\A M)$ and prove  that this homomorphism is either a zero map or  an
isomorphism depending on the choice of a Frobenius covector
(Theorem~\ref{cl-dR}).
   We also show that the complete lift preserves the
Schouten-Nijenhuis bracket of multivector fields
(Proposition~\ref{prop-Sch-br}).
   We prove that for a Frobenius Weil algebra $\A$ it is possible to introduce
a uniquely defined  vertical lift of tensor fields.
   We study the Schouten-Nijenhuis brackets of lifts of multivector fields
(Proposition~\ref{prop-vlbr}).
   In this section  formulas for the complete lift of the tensor product and
relations between the lifts and the Lie derivative are derived.

   In Section~5 we study the complete and the vertical lifts of
the Poisson structure $w$ from a Poisson manifold $(M,w)$ to its Weil bundle
$T^\A M$.
We prove that  the complete and the vertical lifts induce homomorphisms of the
Poisson cohomology spaces
$H_{P}^*(M,w)\to H_{P}^* (T^\A M,w^C)$ and
$H_{P}^*(M,w)\to H_{P}^* (T^\A M,w^V)$
and establish the structure of  these homo\-mor\-phisms for some Poisson
manifolds, in particular, for sym\-plectic manifolds.
   Finally, we compute the modular classes of the lifts of Poisson structures.
   Namely, we prove that the modular class of
Poisson manifold $(T^\A M, w^C)$ is represented by
the vector field $\dim \A \cdot \Delta^V$ for every modular
vector field $\Delta$ of the base manifold $(M,w)$ and that the
modular class of  Poisson manifold
$(T^\A M, w^V)$ is zero (Theorem~\ref{th-mod}).

  The research of this paper  was motivated by the works of
G.~Mitric and I.~Vais\-man~\cite{M-V} and  J.~Grabowski and
P.~Urba\'nski~\cite{G-U, G-U2, G-U3}.

\par

%
%

\section{Poisson manifolds}
\label{sect-pm}

%
%

  Let $M$ be a smooth manifold, $\dim M=m$.
  We will denote by $C^\infty(M)$ the algebra of smooth functions on $M$
and by ${T}^{r,s}(M)$ the space of tensor fields of type~$(r,s)$
on~$M$.
  The algebra of smooth exterior forms on $M$ will be denoted
by~$\Omega^*(M)=\mathop{\oplus}\limits_{k=0}^m \Omega^k(M)$,
and the exterior algebra of skew-symmetric contravariant tensor fields (multivector fields) on~$M$ by~${\cal V}^*(M)=\mathop{\oplus}\limits_{k=0}^m{\cal V}^k(M)$.
  We assume all  manifolds and  maps between
manifolds under consideration to be of class~$C^{\infty}$.
\par
  Throughout the paper we use the Einstein summation convention.
\par
  For  $u\in{\cal V}^k(M)$ we denote by
$i(u):\Omega^p(M)\to  \Omega^{p-k}(M)$
the interior product of a $p$-form with~$u$.
  In local coordinates
\begin{equation}
\label{i(u)}
(i(u)\alpha)_{i_1\dots i_{p-k}}=
u^{j_1\dots j_k}\alpha_{j_1\dots j_ki_{1}\dots i_{p-k}}.
\end{equation}
\par
   We will denote by $d$ the exterior differential on $\Omega^*(M)$
and by~$H^*_{dR}(M)=H^*_{dR}(M,\R)$ the de Rham cohomology of $M$.
\par
  By~${\cal L}_X$ we will denote the Lie derivative along the vector
field~$X$.

  The Lie bracket of vector fields on $M$ can be uniquely extended to
an $\R$-bilinear bracket $[\,\cdot\,, \,\cdot\,]$ on
${\cal V}^*(M)$, called {\it Schouten-Nijenhuis bracket}
\cite{Nij, Schou}, such that $({\cal V}^*(M), [\,\cdot\,, \,\cdot\,])$ is a graded superalgebra.

  The {\it Schouten-Nijenhuis bracket} is an $\R$-bilinear map
$[\,\cdot\,, \,\cdot\,] : {\cal V}^p(M) \times {\cal V}^q(M) \to
{\cal V}^{p+q-1}(M)$ defined as follows.
  Let $X_1$, \dots $X_p$, $Y_1$, \dots, $Y_q$ be vector fields on~$M$.
  Then
$$
\begin{array}{l}
 [X_1\wedge \ldots \wedge X_p,Y_1\wedge \ldots \wedge Y_q] =
  \\[5pt]
\displaystyle    \qquad
    = \sum (-1)^{i+j}  X_1 \wedge \ldots
     \widehat{X_i} \ldots \wedge X_p \wedge [X_i,Y_j]\wedge
     Y_1 \wedge \ldots \widehat{\, Y_j\, } \ldots \wedge Y_q,
\end{array}
$$
where $[X_i,Y_j]$ in the right-hand side
is the Lie bracket of vector fields  and $\widehat{X}$ means the omission of~$X$.

  Let $u\in{\cal V}^p(M)$ and $v\in{\cal V}^q(M)$
be given in terms of a local coordinates as
$u=
u^{i_1\dots i_p} \frac{\partial}{\partial x^{i_1}}
\wedge\ldots \wedge \frac{\partial}{\partial x^{i_p}}$,
$v=
v^{j_1\dots j_q} \frac{\partial}{\partial x^{j_1}}
\wedge\ldots \wedge \frac{\partial}{\partial x^{j_q}}$.
  Then for their Schouten-Nijenhuis bracket we have
\begin{equation}
\begin{array}{rl}
\label{Sch-br}
  [u,v]^{k_2\dots k_{p+q}} = &
\displaystyle
\varepsilon^{k_2\dots k_{p+q}}_{i_2\dots i_p j_1\dots j_q}
u^{ri_2\dots i_p} \frac{\partial}{\partial x^r}
   v^{j_1\dots j_q}\,
+\, (-1)^p
\varepsilon^{k_2\dots k_{p+q}}_{i_1\dots i_p j_2\dots j_q} v^{r
j_2\dots j_q} \frac{\partial}{\partial x^r} u^{i_1\dots i_p}.
\end{array}
\end{equation}
where $\varepsilon^{k_1\dots k_s}_{\ell_1\dots \ell_s}$
is the antisymmetric Kronecker symbol.

  The Schouten-Nijenhuis bracket is supercommutative
$$
  [u,v]=(-1)^{|u|\cdot |v|} \, [v,u],
$$
it satisfies the super-Jacobi identity
\begin{equation}
\label{sJacobi}
  (-1)^{|u| \cdot |v|} \, [[v,y],u] +
  (-1)^{|v| \cdot |y|} \, [[y,u],v] +
  (-1)^{|y| \cdot |u|} \, [[u,v],y] = 0
\end{equation}
and the super-Leibniz identity
\begin{equation}
\label{sLeib}
  [u, v\wedge y] = [u, v]\wedge y  +
   (-1)^{(|u|-1)\cdot |v|}\, v\wedge [u,y],
\end{equation}
where $|u|$ denotes the degree of $u$.
  For more detail see,~e.g.,~\cite{APP, Mich, Nij, Schou, daS-W}.

\smallskip
A {\it Poisson bracket} on a smooth manifold~$M$
is a bilinear skew-symmetric mapping
$\{\, , \,\}:C^\infty(M)\times C^\infty(M)\to C^\infty(M)$,
satisfying the Leibniz rule
$$
  \{f, gh\}=\{f,g\}h+g\{f,h\}
$$
and the Jacobi identity
$$
  \{\{f,g\},h\} + \{\{g,h\},f\} + \{\{h,f\},g\} = 0.
$$
\par
  A  {\it Poisson manifold} is a smooth
manifold $M$ endowed with a Poisson bracket.

  A Poisson bracket on  $M$ uniquely defines a bivector field
$w  \in {\cal V}^2(M)$, called a {\it Poisson bivector},
such that
\begin{equation}
\label{bracket}
\{f,g\}=i(w)(df\wedge dg)
\end{equation}
for any $f,g\in C^\infty(M)$.
  It is known  (see, e.g., \cite{K, Li, Vai-LGPM})
that a bracket on
$C^\infty(M)$ defined by (\ref{bracket}) for a bivector field $w$
satisfies the Jacobi identity if and only if
\begin{equation}
\label{ww0}
[w,w]=0.
\end{equation}
  In terms of local coordinates,~(\ref{ww0}) takes the form
$$
  \displaystyle
  w^{js}\frac{\partial w^{k\ell}}{\partial x^s}+
  w^{ks}\frac{\partial w^{\ell j}}{\partial x^s}+
  w^{\ell s}\frac{\partial w^{jk}}{\partial x^s} =0.
$$
  In what follows  we will denote a Poisson manifold by $(M,w)$.
\par
   A  Poisson bivector determines a bundle map
\begin{equation}
\label{w-map}
  \widetilde w : T^* M \to TM,
\end{equation}
defined by
$$
  \bigl(\widetilde w \alpha\bigr) (\beta) := w(\alpha,\beta),
  \quad \alpha,\beta\in T^*M.
$$
  A Poisson bracket induces a bracket of 1-forms on $M$ by
\begin{equation}
\label{form-bracket}
  \{\alpha,\beta\} = {\cal L}_{\widetilde w \alpha} \beta -
    {\cal L}_{\widetilde w \beta}\alpha - d(w(\alpha,\beta)).
\end{equation}
  This bracket naturally extends the bracket
$\{df, dg\}:= d\{f,g\}$ from
$B^1(M) := \{df \, |\, f\in C^\infty(M)\}$ to~$\Omega^1(M)$, and
$(\Omega^1(M), \{\cdot,\cdot\})$ is a Lie algebra~\cite{K-M,Vai-LP}.
\par
  To each function  $f\in C^\infty(M)$ there is associated a
vector field $X_f=X_f^w\in{\cal V}^1(M)$ called
{\it Hamiltonian vector field of $f$}
 defined by
$X_f(g):=\{f,g\}$.
  Locally, Hamiltonian vector fields on
$(M,w)$ are of the form \cite{M-V, daS-W, Vai-LGPM}
\begin{equation}
\label{Xf}
  X_f^w = \{ f, \cdot \}_w =
     w^{ij}\frac{\partial f}{\partial x^i}
    \frac{\partial}{\partial x^j}, \qquad
      f\in C^\infty(M).  \Bigr.
\end{equation}

  A smooth function $f\in C^\infty(M)$ is called a
{\it Casimir function} if
$\{f,g\}=0$ for every
$g\in C^\infty(M)$, which
is equivalent to the fact that the Hamiltonian vector field
$X_f^w$ of $f$ is trivial.
  In terms of local coordinates, Casimir functions satisfy
$$
   w^{ij} \frac{\partial f}{\partial x^i} = 0.
$$

  {\it The rank\/} of a Poisson bracket  (Poisson structure)
at a point $x\in M$ is defined to be the  rank of $w(x)$.
  The rank of a Poisson bracket on $M$ is the number
$$
   \max\limits_{x\in M} \, \rank  w(x).
$$

 A Poisson manifold~$(M,w)$ is said to be {\it regular}, if
the rank of  $w$ is constant on~$M$.

  A smooth map $\varphi : (M,w) \to (M',w')$  between two Poisson manifolds
is called a {\it Poisson map} if
$$
   \varphi^*(\{f,g\}_{M'})=\{\varphi^*(f),\varphi^*(g)\}_M.
$$
  There is also an alternative characterization of Poisson maps~\cite{daS-W}.
  Let $X\in {\cal V}^k(M)$ and $Y\in{\cal V}^k(M')$ be two multivector fields.
  We say that $X$ is {\it $\varphi$-related\/} to $Y$, writing $Y=\varphi_* X$, if
$$
 (\wedge^k T_x \varphi) X(x) = Y(\varphi(x)) \quad \hbox{for all } x\in M.
$$
  Then $\varphi$ is a Poisson map if and only if
\begin{equation}
\label{pmap}
   w'=\varphi_* w.
\end{equation}

\medskip
 For a Poisson manifold $(M,w)$ A.\,Lichnerowicz~\cite{Li} introduced an operator
$$
\sigma = \sigma_w : {\cal V}^k(M) \to {\cal V}^{k+1}(M)
$$
defined by $\sigma u := [w,u]$.
  From the super-Jacobi identity (\ref{sJacobi}) and the super-Leibniz rule
(\ref{sLeib}) it follows that
$$
  \sigma( u\wedge v) =
    \sigma u\wedge v + (-1)^{|u|} u \wedge \sigma v
$$
and
$$
\sigma \circ \sigma =0.
$$
Therefore the {\it Poisson cohomology spaces}
$$
\displaystyle
  H^k_{P} (M,w) := \frac{\ker \sigma: {\cal V}^k(M) \to
   {\cal V}^{k+1}(M)}
   {\im \sigma: {\cal V}^{k-1}(M) \to {\cal V}^{k}(M)}.
$$
are defined.
  In the general situation this cohomology is very difficult to compute (see, e.g.,
\cite{Gam, G-L, Mon, Naka, Pich, R-V, Vai-LP, Xu}).
\par
   A map  $\widetilde w$ can be extended to a map
$\Omega^k(M)\to {\cal V}^k(M)$ defined by the formula
\begin{equation}
\label{tw}
   \widetilde w \theta(\alpha_1, \dots, \alpha_k) =
   (-1)^k \theta (\widetilde w\alpha_1, \dots,
     \widetilde w\alpha_k),
\end{equation}
where $\alpha_i \in \Omega^1(M)$.
   In terms of local coordinates,
\begin{equation}
\label{tw-f}
  (\widetilde w \theta)^{j_1\dots j_k}=
(-1)^k w^{i_1j_1}\dots w^{i_kj_k} \theta_{i_1\dots i_k}.
\end{equation}
  Clearly, for a symplectic manifold map~(\ref{tw})
is an isomorphism.
  It can be shown that \cite{Vai-LP}
\begin{equation}
\label{d-sigma}
\sigma \circ \widetilde w = (-1)^k \widetilde w \circ d.
\end{equation}
   It follows that there arise natural homomorphisms
$$
\rho^k : H^k_{dR}(M) \to H^k_{P}(M,w).
$$
  In the case of a symplectic manifold, these homomorphisms are isomorphisms,
and the Poisson cohomology is isomophic to the de Rham cohomology.
  See~\cite{K, Li, Vai-LP} for details.

\begin{example}
  Let $M$ be a smooth manifold with  zero Poisson structure ($w=0$).
  Then
$$
   H^k_P(M,w) \cong {\cal V}^k(M).
$$
\end{example}

\begin{example}
\label{ex-reg-pr}
  Let $S$ be a symplectic  manifold and $N$  an arbitrary smooth
manifold.
  Let  $M=S\times N$ be the regular Poisson manifold whose Poisson structure~$w$
is induced  from~$S$.
  Suppose that $\dim H^*_{dR}(S) < \infty$.
  Then~\cite{Vai-LP, Vai-LGPM}
\begin{equation}
\label{reg-pr}
H_P^r(M,w) \cong \mathop{\oplus}\limits_{0\le k\le r} H^k_{dR}(S)
\otimes {\cal V}^{r-k}(N).
\end{equation}
\end{example}

\par
   For a Poisson manifold $(M,w)$, the
{\it canonical cohomology class}
$[w]\in H^2_P(M,w)$ is defined.
  This class is zero if and only if there exists
$X\in {\cal V}^1(M)$ such that ${\cal L}_X w =w$.

  A vector field  $X$ such that  ${\cal L}_X w =w$ is called a
{\it Liouville vector field for $w$}~\cite{daS-W}.
  A Poisson manifold $(M,w)$ admitting a  Liouville vector field is said to be
{\it exact} (or {\it homogeneous})~\cite{Vai-LGPM}.

\par
  Assume that $M$ is oriented, and let  $\mu$
be a volume form on $M$.
  The divergence  ${\rm div}_\mu X$
of a vector field $X\in {\cal V}^1(M)$ is the smooth function defined by
$$
  {\cal L}_X \mu = ({\rm div}_\mu X)\mu.
$$
  For a Poisson manifold $(M,w)$ with  volume form~$\mu$,
the operator
$$
  \Delta_\mu = \Delta_{\mu, w} : f\in C^\infty(M) \, \longmapsto\,
   {\rm div}_\mu X_f\in C^\infty(M)
$$
is a derivation on $C^\infty(M)$, hence a vector field~\cite{Wein-m}.
  This vector field is called the {\it modular vector field\/}
of $(M,w,\mu)$.
\par
  The modular vector field satisfies $\sigma
\Delta_\mu=0$~\cite{KS}.
  For another  volume form
$a\mu$, where $a\in C^\infty(M)$ is a non-vanishing function,
the modular vector field changes to
$\Delta_{a\mu} =
\Delta_\mu + X_{-{\rm log}\, a}$~\cite{Wein-m}.
  Since Hamiltonian vector fields form the space of
1-coboundaries of~$\sigma$~\cite{daS-W, Li}, it follows that the set
of modular vector fields for all
volume forms on $M$ is a cohomology class from $H_{P}^1(M,w)$.
  This cohomology class is called the {\it modular class\/}
of~$(M,w)$.
  We denote it by ${\rm mod}(M,w)$.

  From~(\ref{Xf}) it follows that
in terms of local coordinates  on $M$
the modular vector field is given by~\cite{M-V}
\begin{equation}
\label{mod-vf}
  \Delta_{\mu} =
   \sum_{j=1}^m \left(
    \frac{\partial w^{ij}}{\partial x^j} +
    w^{ij}  \frac{\partial\, \log \rho}
        {\partial x^j} \right) \frac{\partial }{\partial x^i},
\end{equation}
where $\mu=\rho\, dx^1\wedge \ldots \wedge dx^m$.

  In the case when $M$ is non-orientable, one defines the modular class
in a similar way using a smooth density instead of a volume form.

%
%

\section{Weil algebras and Weil bundles}
\label{sect-wa}

\subsection{Weil bundles as smooth manifolds over Weil algebras}
\label{subsec-wb}

\par
  A {\it Weil algebra} \cite{KMS, VSh-ljm} is a finite-dimensional
associative commutative $\R$-algebra
$\A$ with unit $1_\A$
whose nilpotent elements form a unique maximal ideal
$\Ao$.
  The linear span of $1_\A$ form a subalgebra isomorphic to $\R$.
  We will identify it with $\R$.
  As a vector space, $\A$ is the direct sum $\R\oplus{\Ao}$.
  In what follows  $n=\dim_\R \Ao$ and so $\dim \A=n+1$.
\par
  By $\Ao{}^r$ we denote the $r$th power of~$\Ao$.
  The positive integer $h$ defined by the relations
$\Ao{}^h\ne 0$, $\Ao{}^{h+1}=0$
is called the  {\it height} of $\A$.
  Let $d_k(\A)={\rm dim}_\R\,{{\Ao{}^k}/{\Ao}{}^{k+1}}$ for $k=1,\dots, h$ and
$d_0(\A)={\rm dim}_\R \,\A/\Ao =1$.
  The number $d_1(\A)$ is usually called the
{\it width} of $\A$.
\par
  The chain of embedded ideals
$$
\A \supset {\Ao} \supset \Ao{}^2 \supset \ldots \supset
\Ao{}^h \supset 0
$$
can be extended to the chain of ideals called
the Jordan-H\"older composition series~\cite{VSh-ljm}
$$
\A \supset {\Ao} = \I_1 \supset \I_2 \supset \ldots \supset
\I_n \supset 0,
$$
where $\I_a/\I_{a+1}$ is a 1-dimensional algebra with  zero
multiplication.
  Here
$$
\Ao{}^k=\I_{1+d_1(\A)+\dots+d_{k-1}(\A)}\quad \hbox{for} \quad 2\le k\le h.
$$
  This is a particular case of the general construction for  rings,
see~\cite{Pier}.
  Using the Jordan-H\"older composition series one can choose a
basis  {\it (a Jordan-H\"older basis)}
\begin{equation}
\label{jgb}
\{e_a\}=\{e_0,e_{\hat a}\}, \quad
a=0,1,\dots,n = \dim \Ao,\quad \hat a =1, \dots, n,
\end{equation}
in $\A$ such that
 $e_0=1_\A \in \R$, $e_{\hat a}\in {\I}_{\hat a}$,
$e_{\hat a}\not \in{\I}_{\hat a +1}$.
  For $X=x^ae_a=x^0+x^{\hat a}e_{\hat a}\in \A$ we set
$\Xo=x^{\hat a}e_{\hat a}$, then $X=x^0+\Xo$.
  Let $\delta^a$ be the coordinates of unit of $\A$, i.e.,
$1_\A=\delta^a e_a$.
  Sometimes we  denote the multiplication in $\A$ by a dot in order
to avoid confusion.
\par
  We denote by $(\gamma_{ab}^c)$
the structure tensor of  $\A$ with respect to a Jordan-H\"older basis~(\ref{jgb}).
  We have
$e_ae_b= \gamma_{ab}^c e_c$,
$\gamma^b_{0a}=\delta^b_a$ (Kronecker's deltas), and
$\gamma^c_{a\hat b}=0$ for $a\ge c$.
  The conditions of commutativity and associativity are
$\gamma_{ab}^c=\gamma_{ba}^c$ and
\begin{equation}
\label{ass}
\gamma_{ab}^c\gamma_{ef}^b=\gamma_{ae}^b\gamma_{bf}^c,
\end{equation}
respectively.
  The relations $e_a=e_a\cdot 1_{\A} = e_a \delta^b
e_b = \delta^b \gamma^c_{ab} e_c$ imply that
\begin{equation}
\label{gamma-delta}
 \delta^b \gamma^c_{ab} = \delta^c_a.
\end{equation}
\par
  A smooth function $f:U\subset\A\to\A$ is said to be
{\it $\A$-differentiable ($\A$-smooth)} if its differential $df$
is an
$\A$-linear map.
  The conditions of $\A$-differentiability of $f$,
usually called  {\it Scheffers' equations}, are (see~\cite{Schef, Shir66, VSS}):
\begin{equation}
\frac{\partial f^{b}}{\partial x^c} \gamma^c_{ad}=
\gamma^{b}_{ac} \frac{\partial f^{c}}{\partial x^d}.
\label{u-s1}
\end{equation}
  Scheffers' equations are equivalent to
\begin{equation}
\frac{\partial f^{b}}{\partial x^a} =
\gamma^{b}_{ac}\delta^d \frac{\partial f^{c}}{\partial x^d}.
\label{u-s2}
\end{equation}
  Let $\A^m=\A\times \cdots\times \A$ be the  $\A$-module
of $m$-tuples of elements of  $\A$.
  We will enumerate the real coordinates in $\A^m$ corresponding
to a basis~(\ref{jgb}) by the double indices~${ia}$.
   For a smooth function of several variables $f: U \subset \A^m \to \A$,
$f:\{X^i=x^{ia}e_a\}\mapsto f(X^i)=f^b(x^{ia})e_b$, Scheffers' conditions
of $\A$-differentiability are of the form~\cite{VSh-ljm,VSS}:
\begin{equation}
\frac{\partial f^{b}}{\partial x^{ia}} =
\gamma^{b}_{ac}\delta^d\frac{\partial f^{c}}{\partial x^{id}}.
\label{u-s3}
\end{equation}
  If  $f$ satisfies (\ref{u-s3}), its differential can be
represented in the form  $df=f_idX^i$,
where $f_i=\delta^a\frac{\partial f}{\partial x^{ia}}$ is the  partial derivative with respect to the variable~$X^i\in\A$.
  We will denote the latter by
$\frac{\partial f}{\partial X^i}$.
  Thus,
\begin{equation}
\label{der-dir}
f_i=\frac{\partial f}{\partial X^i} =
\delta^a\frac{\partial f}{\partial x^{ia}}.
\end{equation}
  The functions $f_i(X^j)$, $i=1,\dots,m$, are also
$\A$-differentiable.
\par
  Recall that a smooth map $f:M\to N$ of a foliated manifold
$(M,{\cal F})$ is called {\it projectable} (or {\it basic}) if
$f$ is constant along the leaves of  $\cal F$.
  The natural epimorphism $\pi^m:\A^m\to\R^m$ determines
{\it the canonical $\Ao{}^m$-foliation} on $\A^m$.
  The following theorem (see~\cite{VSh-ljm}) describes the local
structure of an  $\A$-differentiable map of the form
$F:U\subset\A^m\to\A^k$ for a Weil algebra $\A$.
\par
\begin{theorem}
{\bf~\cite{VSh-ljm}~}
{\rm 1)}
Let $U\subset\A^m$ be an open set
and $\varphi:U\to\A^k$ a projectable map with respect to
the canonical $\Ao{}^m$-foliation. Then the formula
\begin{equation}
X^{i'}=\varphi^{i'}+\sum^h_{|p|=1}\frac{1}{p!}
\frac{D^p\varphi^{i'}}{Dx^p}\Xo{}^p,
\label{A-prol}
\end{equation}
where $i=1,\dots,m$,  $i'=1,\dots,k$,  $p=(p_1,\dots,p_m)$ is a multiindex of length~$m$ and $p!=p_1!\dots p_m!$, $X^i=x^i+\Xo{}^i$ is the decomposition with respect to {\rm (\ref{jgb})},
$\Xo{}^p=(\Xo{}^1)^{p_1}\dots (\Xo{}^m)^{p_m}$,
determines an $\A$-smooth map
$\Phi:U\to\A^k$.
\par
{\rm 2)} Any projectable $\A$-differentiable map
$\Phi: U\to  \A^k$
is of the form {\rm (\ref{A-prol})} for some basic functions
$\varphi^{i^{\prime}}:U\to\A$.
\end{theorem}
\begin{definition}
 Let  $\varphi:U\to\A^k$ be a projectable map.
Then the map $\Phi:U\to\A^k$ given by
(\ref{A-prol}) is called  the {\it analytic prolongation}
($\A$-{\it prolon\-gation}) of $\varphi$.
\end{definition}

 The analytic prolongation of a map
$\varphi$ will be denoted by $\varphi^\A$.

\begin{proposition}
\label{A-prol-prop}
{\bf \cite{VSh-ljm}~}
Analytic prolongations satisfy the following relations:
\par
$1^{\circ}$.
$(\varphi+\psi)^{\A}=\varphi^\A+\psi^\A$.
\par
$2^{\circ}$.
$(\varphi\cdot\psi)^{\A}=\varphi^\A\cdot\psi^\A$.
\par
$3^{\circ}$.
$(\varphi^\A\circ\psi)^{\A}=\varphi^\A\circ\psi^\A$.
\par
$4^{\circ}$.
$({D^p\varphi}/{Dx^p})^\A= {D^p\varphi^\A}/{DX^p}$~
for~ $\varphi: U\subset \A^m\to \A$.
\end{proposition}

  Let now  $\Lm$ be an arbitrary  $\A$-module, $\dim_\R\Lm<\infty$,
and let $M$ be a smooth manifold.
\par
  An $\Lm$-{\it chart} on $M$ is a pair $(U, h)$ consisting of an open set
$U\subset M$ and a diffeomorphism
$h:U\to U'\subset\Lm$.
  An $\Lm$-{\it atlas} on $M$ is a collection of $\Lm$-charts
$\{(U_{\kappa}, h_{\kappa})\}_{\kappa\in K}$ such that
$\{U_{\kappa}\}_{\kappa\in K}$ is a covering of~$M$ and the tangent map
\begin{equation}
\label{tan-iso}
T_{h_{\lambda}(x)}(h_{\kappa}\circ h_{\lambda}^{-1}):
T_{h_{\lambda}(x)}\Lm\equiv\Lm\to\Lm\equiv T_{h_{\kappa}(x)}\Lm
\end{equation}
is an isomorphism of $\A$-modules for all
$x\in M$ and $\kappa,\lambda \in K$.
  This condition is equivalent to the $\A$-differentiability of all
transition functions $h_{\kappa\lambda}:=h_{\kappa}\circ
h_{\lambda}^{-1}$.

\begin{definition}
{\bf\cite{VSS}}
  A quadruple $(\A,\Lm,M,{\cal A})$ consisting of a Weil algebra $\A$,
an $\A$-module $\Lm$, a smooth manifold $M$, and a maximal
$\Lm$-atlas $\cal A$ on~$M$ is called an {\it $\Lm$-manifold} or an {\it
$\A$-smooth manifold, modeled on\/ $\Lm$}.
\end{definition}
\par
   The isomorphism (\ref{tan-iso}) allows ones to transfer by means of
$T_x h_{\kappa}^{-1}$ the structure of an $\A$-module from $T_x U'_{\kappa}\cong \Lm$
to the tangent space $T_xM$ at $x\in M$.

\par
\begin{definition}
{\bf\cite{VSS}}
  $\A^n$-manifold is called an {\it $n$-dimensional $\A$-smooth
manifold}.
\end{definition}
\par
  Let  $M$ and $M'$ be two $\A$-smooth manifolds modeled, respectively, by
$\A$-modules $\Lm$ and $\Lm'$.
  A smooth map $f:M\to M'$ is said to be {\it $\A$-smooth}
if $T_xf$ is an $\A$-linear map for all
$x\in M$.
  This is equivalent to the $\A$-differentiability
of all chart representatives $h_{\kappa'}\circ f\circ h_\kappa^{-1}$
of~$f$.


  For any Weil algebra $\A$ and smooth manifold $M$, the {\it Weil
bundle\/} $\pi_\A:T^\A M \to M$ of $\A$-points is defined as
follows.
  An $\A$-point near to  $x\in M$ is a homomorphism
$X:C^\infty(M)\to \A$ such that the real part of $X(f)\in\A$
coincides with $f(x)$.
    The set $T^\A M$ of all $\A$-points near to points of $M$ can be
endowed with a structure of smooth manifold.
  Bundle projection $\pi_A$ sends $X$ to $x$.
  Let $x^i$ be local coordinates on a neighborhood
$U$ of $M$.
  These coordinates induce local coordinates $x^{ia}$ on $\pi_\A^{-1}(U)\subset T^\A M$
defined by
$X(x^i)=x^{ia}(X) e_a$, $a=0,1,\dots n=\dim \Ao$, where $x^{i0}=x^i$.
  The functions $X^i=x^{ia} e_a$ are $\A$-valued  coordinates
on $\pi_\A^{-1}(U)$.
  A.P.~Shirokov proved~\cite{Shir12} that these coordinates define the structure
of a smooth manifold over $\A$ on  $T^\A M$.

  The correspondence which assigns to a manifold $M$ the Weil
bundle $T^\A M$ and to a smooth map $\varphi:M\to N$ the map
$\varphi^\A: T^\A M \ni X \mapsto X\circ \varphi^* \in T^\A N$,
where $\varphi^*: C^\infty (N) \ni g\mapsto g \circ \varphi\in
C^\infty(M)$, is a functor called the {\it Weil functor}~
(see~\cite{KMS, VSh-ljm, Weil}).
  It is well-known that Weil functors preserve products, i.e., $T^\A(M\times N)\cong
T^\A M\times T^\A N$.


\subsection{The structure of a Frobenius Weil algebra}
\label{subsec-struc-fwa}

\par
\begin{definition}
{\bf\cite{Cur-R}}~A {\it Frobenius Weil algebra\/} is a pair
$(\A,q)$ where $\A$ is a Weil algebra and
$q:\A\times\A\to \R$ is a nondegenerate bilinear form which satisfies the following associativity condition:
\begin{equation}
\label{q1}
q(XY,Z)=q(X,YZ) \quad
  \hbox{for any} \quad X,Y,Z\in \A.
\end{equation}
\end{definition}
\par
  In terms of basis (\ref{jgb}) the condition (\ref{q1})
is written as
\begin{equation}
\label{q2}
q_{ac}\gamma_{bd}^c=\gamma_{ab}^c q_{cd}.
\end{equation}

  The form $q$ is called the {\it Frobenius form}.
  Frobenius algebras play an important role in the theory of smooth
manifolds over algebras in constructing realizations of tensor
operations \cite{VSS, Kruch}.

  For a Frobenius algebra $\A$ the linear form
$p:\A\to\R$ defined by
\begin{equation}
\label{p-def}
p(X)=q(X,1_\A)
\end{equation}
is called the {\it Frobenius covector}.
  Its coordinates  satisfy
\begin{equation}
\label{p}
p_a\gamma_{bc}^a=q_{bc}.
\end{equation}
  Contracting (\ref{p}) with~$\delta^c$  gives
\begin{equation}
\label{p2}
p_b = q_{bc} \delta^c.
\end{equation}
  From (\ref{q1}) and (\ref{p}) it  follows that
\begin{equation}
\label{p-q}
q(X,Y)=p(XY) \quad  \hbox{for any} \quad X,Y\in \A.
\end{equation}
\par

\par
   Let $\A^*$ be the dual space of $\A$, i.e., the space of $\R$-linear functions $\xi:\A\to\R$.
  The Frobenius form $q$ induces the isomorphism
\begin{equation}
\label{phi}
\varphi: \A\to \A^*, \qquad
\varphi(X)Y := q(X,Y) = p(XY).
\end{equation}
  From the definition of $p$ it follows that
$\varphi(1_\A)=p$, and (\ref{q1})~implies that
$$
  \varphi(XY)(Z) = \varphi(X)(YZ).
$$

  The isomorphism $\varphi$ allows ones to transfer
the multiplication operation from $\A$ to~$\A^*$:
$$
  \xi * \eta := \varphi(\varphi^{-1}(\xi) \cdot
      \varphi^{-1}(\eta))
$$
(the dot means the multiplication in~$\A$).
  The multiplication $*$ turns $\A^*$ into a Weil algebra isomorphic to~$\A$.
    The bilinear form
$q$ induces the form $\widetilde q : \A^* \times \A^* \to \R$,
$$
  \qw (\xi, \eta) := \xi (\varphi^{-1}\eta).
$$

It is obvious that
   $\qw$ is symmetric.

   Let $\{\ea^a\}$ be the basis in $\A^*$ dual to
a basis~$\{e_a\}$ of $\A$.
   Then $\varphi(e_a) = q_{ab} \ea^b$, $\varphi^{-1}(\ea^a)
= \qw^{ab} e_b$ and~$\|\qw^{ab}\| =
\|q_{ab}\|^{-1}$.
  In what follows we will omit the tilde in $\qw^{ab}$ and write  $q^{ab}$ instead
of $\qw^{ab}$.

\begin{proposition}
  For any $\xi,\eta\in\A^*$
\begin{equation}
\label{q-inv-1}
   \qw(\xi, \eta) = (\xi * \eta)(1_\A).
\end{equation}
\end{proposition}

\begin{Proof}
  Let $\xi = \varphi(X)$, $\eta = \varphi(Y)$ for some
$X,Y\in  \A$.
  Then
$\qw(\xi, \eta) = \xi(\varphi^{-1}\eta) = \varphi(X) Y =
p(XY) = p (\varphi^{-1} (\xi) \cdot\varphi^{-1}(\eta)) =
p(\varphi^{-1}(\xi*\eta))$.
  Let us show that $p(\varphi^{-1}\zeta) = \zeta(1_\A)$ for any
$\zeta\in\A^*$.
  In fact, if $\zeta = \varphi(Z)$, then
$p(\varphi^{-1}\zeta) = p(Z) = \varphi(Z)(1_\A) = \zeta(1_\A)$.
\end{Proof}

   Every function $F:\A \to \A$ obviously determines the function
$$
\Fs = \varphi \circ F \circ \varphi^{-1} : \A^* \to \A^*.
$$

  The proof of the following proposition is immediate.

\begin{proposition}
   If $F$ is $\A$-differentiable, then $\Fs$ is $\A^*$-differentiable.
\end{proposition}

\par
\begin{example}
  Let $\A$ be the Weil algebra of plural numbers
$$
{\D}^n=\R(\varepsilon^n)=\{x_0+x_1\varepsilon+\dots+x_{n}\varepsilon^{n}
\,| \, x_i\in{\R}, \varepsilon^{n+1}=0\}
$$
(the algebra of truncated polynomials in one variable
$\varepsilon$ of degree not greater than~$n$).
  For  $n=1$ we get the algebra of dual numbers (also known as Study numbers)
$$
{\D}=\R(\varepsilon)=\{x_0+x_1\varepsilon\,|\, x_0,x_1\in{\R},
\varepsilon^2=0\}.
$$

  There is a natural Jordan-H\"older basis in $\R(\varepsilon^n)$,
namely, $e_0=1$, $e_a=\varepsilon^a$, $a=1,\dots,n$.
  Let
$p=(p_0,\dots, p_n)$ be an arbitrary covector on
$\R(\varepsilon^n)$.
  Then the matrix $\|p_c \gamma_{ab}^c\|$  is
\begin{equation}
\label{pl-fr}
\left(
\begin{array}{lllll}
p_0 & p_1 &\dots & p_{n-1} & p_n\\
p_1 & p_2 & \dots & p_n & 0\\
\,\vdots & \,\vdots & \backddots & \backddots & \,\vdots\\
p_{n-1} & p_n & \backddots & 0 & 0\\
p_n & 0 & \dots & 0 & 0  \\
\end{array}
\right).
\end{equation}

  Since $\det\|p_c\gamma_{ab}^c\|=p_n^{n+1}$, it follows that
$p$
is a Frobenius covector if and only if
$p_n=p(e_n)\ne 0$ (cf. Proposition~\ref{prop-fa}).

\end{example}

\par
  In what follows we assume all Weil algebras under consideration to
be Frobenius algebras.


  Let $(\A,q)$ be a Frobenius Weil algebra of height $h$ and let $p$ be
its Frobenius covector.
  We denote $n=\dim \Ao$.
  Let us fix a Jordan-H\"older basis (\ref{jgb}) in $\A$.

\begin{proposition}
\label{prop-fa}
  For a Frobenius Weil algebra $(\A,q)$ the following conditions hold:

  {\rm 1)} $\dim\Ao{}^h=1$, that is, $\Ao{}^h=\I_n${\rm ;}

  {\rm 2)} $p|_{\Ao{}^h}\ne 0$.
\end{proposition}
\begin{Proof}
  Denote
$$
\Ann\Ao:=\{X\in \Ao \,|\, X \cdot \Ao =0\}.
$$
  Let $0\ne X\in\Ann\Ao$.
  Then for any $Y=y^0+\Yo\in \A$ we have
$X Y = X y^0$ and $q(X,Y)= p(X Y)= y^0 p(X)$.
  The nondegeneracy of $q$ implies that $p(X)\ne0$.
  Hence $\Ann \Ao \cap \ker p=0$, which means that
$\dim\Ann \Ao\le 1$.
  However, it is obvious that $0\ne \Ao{}^h\subset \Ann \Ao$.
  Consequently,
$\dim \Ao{}^h=\dim\Ann \Ao= 1$.
\end{Proof}

  The second statement of Proposition~\ref{prop-fa} is equivalent to the
inequality $p_n\ne 0$.
  In what follows we will always  assume that a Jordan-H\"older basis~(\ref{jgb})
is chosen in such a way that $p_n=p(e_n)=1$.
  Then the matrix of~$q$ is of the form
$$
\|q_{ab}\| =
\left(
\begin{array}{ccccc}
*      & * & \dots & * & 1\\ \relax
*      & * & \dots & * & 0\\
\vdots & \vdots & \ddots & \vdots & \vdots\\ \relax
*      & * & \dots & * & 0\\
1      & 0 & \dots & 0 & 0
\end{array}
\right),
$$
and the inverse matrix~$\|q^{ab}\|$ is of the form
\begin{equation}
\label{q-inv}
\|q^{ab}\| =
\left(
\begin{array}{ccccc}
0      & 0 & \dots & 0 & 1\\ \relax
0      & * & \dots & * & *\\
\vdots & \vdots & \ddots & \vdots & \vdots\\ \relax
0      & * & \dots & * & *\\
1      & * & \dots & * & *
\end{array}
\right).
\end{equation}

\par
\begin{remark}
\label{rem2}
  Let $\po\in\A^*$ be defined by
$\po(X):=x^0$ (i.e., the projection onto $\R$ along $\Ao$).

 We show that if  $p\in\A^*$ is a Frobenius covector
then  $\widetilde p:=p-p(1)\po$
is also a Frobenius covector.
  Suppose the contrary.
  Then there exists $X\in \A$ such that $\widetilde p(XY)=0$
for any $Y\in \A$.
  This means that $x^0p(\Yo) + y^0p(\Xo) + p(\Xo\Yo)=0$.
  Let  $Z\in \Ao{}^h$, $p(Z)=1$
(in terms of a Jordan-H\"older basis under consideration, $Z=e_n$)
and let
$\widetilde X=X-x^0 Z$.
  Then for any $Y\in \A$ we have
$\widetilde X Y = XY - x^0 y^0 Z = x^0 y^0 + x^0\Yo + y^0\Xo + \Xo\Yo
- x^0 y^0 Z$.
  Obviously, $p(\widetilde X Y)=0$, which
contradicts to the fact that $p$ is a Frobenius covector.
  This means that having an arbitrary Frobenius covector $p\in \A^*$ we can obtain another
Frobenius covector $\widetilde p$ satisfying
$\R \subset \ker \widetilde p$.
\end{remark}

  Let $\{e^a\}$ be the basis in $\A$, corresponding to the
dual basis $\{\ea^a\}$  in $\A^*$ with respect to the isomorphism
$\varphi$.
Then
\begin{equation}
\label{jgb-dual}
e^a=q^{ab} e_b, \quad e_a=q_{ab} e^b.
\end{equation}
  Denote by~$\gamma^{ab}_c$  the  components of the structure tensor of~$\A$ with
respect to the basis~$\{e^a\}$.
  Then
\begin{equation}
\label{pee-up}
 p(e^a e^b)=q^{ab}.
\end{equation}

  The components $\gamma_{ab}^c$ and $\gamma^{ab}_c$ are related
as follows
$\gamma^{ab}_c= q^{ad} q^{bf} q_{ch} \gamma^h_{df}$.
  Using~(\ref{q2}), we  have
\begin{equation}
\label{gamma-rel}
  \gamma^{ab}_c = q^{ad}  \gamma^b_{dc},
\quad
  \gamma^b_{dc} =  \gamma^{ab}_c  q_{ad}.
\end{equation}
From (\ref{gamma-delta}) it follows that
\begin{equation}
\label{gamma2-delta}
  \gamma^{ab}_c \delta^c = q^{ab}.
\end{equation}
  From (\ref{gamma-rel}) and (\ref{ass}), we have
\begin{equation}
\label{gamma-gamma}
 \gamma^{af}_c \gamma^d_{bf} = \gamma^{ad}_f \gamma_{bc}^f.
\end{equation}
   Using
(\ref{gamma-rel}), one can obtain the following formulas for the
products $e_a e^c$:
\begin{equation}
\label{ee-gamma}
 e_a e^c = \gamma_{ab}^c e^b = \gamma^{cb}_a e_b.
\end{equation}
  In fact, we have $e_a e^c = q_{ad} e^d e^c =
q_{ad}\gamma^{cd}_b e^b = \gamma_{ab}^c e^b$. The second relation
is derived in the similar way. From (\ref{ee-gamma}) it follows
that
\begin{equation}
\label{pee}
  p(e_a e^c) =\delta_a^c
\end{equation}
and
\begin{equation}
\label{p-gamma}
  p(e_a e_b e^c) =\gamma_{ab}^c, \quad p(e^a e^b e_c)= \gamma^{ab}_c.
\end{equation}

  Contracting (\ref{p2}) with $q^{ab}$ we obtain
\begin{equation}
\label{q-p}
  q^{ab} p_b =\delta^a.
\end{equation}

   The following two formulas~(\ref{F-contr}) and~(\ref{partialF})
for an $\A$-smooth function
$F:\A\to \A$, $F=F^a(x^b) e_a = F_a(x^b) e^a$ will be used in the
sequel.
\begin{equation}
\label{F-contr}
 F^a p_a = F_b \delta^b.
\end{equation}
\begin{equation}
\label{partialF}
\frac{\partial (\delta^a F_a)}{\partial x^b} =
\delta^c  \frac{\partial F_b}{\partial x^c}.
\end{equation}
  The proof of~(\ref{F-contr}) is obvious.
  Formula~(\ref{partialF}) follows from~(\ref{u-s2}) and~(\ref{gamma-delta}):
$\frac{\partial (\delta^a F_a)}{\partial x^b} =
\delta^a \delta^c \gamma^d_{ab}
\frac{\partial F_d}{\partial x^c} =
\delta^c \frac{\partial (\delta^d_b F_d)}{\partial x^c} =
\delta^c \frac{\partial F_b}{\partial x^c}$.

%
%

\section{Lifts of tensor fields to Weil bundles}
\label{sect-lifts}

\par
  Let $(\A,q)$ be a Frobenius Weil algebra of height~$h$, and
let $p$ be the corresponding Frobenius covector on~$\A$.
    We will use the indices $i$, $j$, \dots, and $\alpha$, $\beta$,\dots{} to enumerate
coordinates on manifolds and the indices $a$, $b$, $c$, $d$,\dots to
enumerate coordinates in~$\A$.
\par
    The Weil bundle $T^\A M$ of an  $m$-dimensional smooth manifold $M$
is an $m$-dimen\-si\-onal $\A$-smooth manifold.
  For
each $X\in T^\A M$ the tangent space  $T_X T^\A M$ is an
$m$-dimensional $\A$-module.
  This allows us to consider $\A$-linear $\A$-valued tensors at any point
$X\in T^\A M$ and $\A$-smooth tensor fields on $T^\A M$ (see, e.g.,~\cite{VSh93}).
\par
   For an arbitrary  $\A$-smooth manifold $M_\A$ the following notation is
used:

  ${T}^{r,s}M_\A$ is the bundle of  $\R$-valued tensors of type~$(r,s)$,

  $\A\otimes {T}^{r,s}M_\A$ is the bundle of $\A$-valued tensors of type~$(r,s)$,

  ${T}_{\Alin}^{r,s}M_\A$ is the subbundle in
$\A\otimes {T}^{r,s}M_\A$ consisting of $\A$-linear tensors.

  Local sections of these bundles (local $(r,s)$-tensor fields on $M_\A$)
are expressed in terms of local coordinates $X^i=x^{ia} e_a$ on
$M_\A$ as follows
$$
  t_{i_1a_1\dots\, i_ra_r}^{j_1b_1\dots j_sb_s} dx^{i_1a_1}\otimes_\R
  \ldots \otimes_\R dx^{i_ra_r} \otimes_\R
   \frac{\partial}{\partial x^{j_1b_1}} \otimes_\R  \ldots
   \otimes_\R \frac{\partial}{\partial x^{j_sb_s}},
$$
where in the first case
$t_{i_1a_1\dots\, i_ra_r}^{j_1b_1\dots j_sb_s}$
are smooth  real-valued functions of $x^{kc}$ , and in the second
and the third cases they are smooth $\A$-valued functions
$t_{i_1a_1\dots\, i_ra_r}^{j_1b_1\dots j_sb_s}=
{\stackrel{c}{t}{}_{i_1a_1\dots\, i_ra_r}^{j_1b_1\dots j_sb_s}} e_c$.
  Local sections of the third bundle
can also be expressed in the form
\begin{equation}
\label{A-ten}
  t_{i_1\dots i_r}^{j_1\dots j_s} dX^{i_1}\otimes_\A
  \ldots \otimes_\A dX^{i_r} \otimes_\A
   \frac{\partial}{\partial X^{j_1}} \otimes_\A  \ldots
   \otimes_\A \frac{\partial}{\partial X^{j_s}}.
\end{equation}
\par
   The total space of the bundle ${T}_{\Alin}^{r,s}M_\A$
carries the natural structure of an $\A$-smooth manifold.
\begin{definition}
   $\A$-{\it smooth tensor field of type $(r,s)$} on $M_{\A}$ is an
$\A$-smooth section of the bundle ${T}_{\Alin}^{r,s}M_\A$.
\end{definition}
   For an $\A$-smooth tensor field functions $t_{i_1\dots i_r}^{j_1\dots j_s}$
are $\A$-smooth functions of $X^k$.

   We denote the space of  $\A$-smooth tensor fields of type $(r,s)$ on~$M_\A$ by
${T}^{r,s}_{\Adiff}(M_\A)$,
the space of $\A$-smooth exterior forms on~$M_\A$
by $\Omega^*_{\Adiff}(M_\A)$ and the space of
$\A$-smooth multivector fields
by~${\cal V}^*_{\Adiff}(M_\A)$.

   For an $\A$-smooth exterior form $\Theta=\Theta_{i_1\dots i_k}
dX^{i_1}\wedge\ldots\wedge dX^{i_k}\in \Omega^*_{\Adiff}(M_\A)$
and $\A$-smooth multivector fields $U\in {\cal
V}_{\Adiff}^g(M_\A)$,
$V\in{\cal V}^\ell_{\Adiff}(M_\A)$  the
exterior differential $d\Theta$ and the Schouten-Nijenhuis bracket
$[U,V]$
can be represented in terms of coordinates $X^i$, respectively, as
follows.
$$
d\,\Xi=\frac{\partial \Xi_{i_1\dots i_k}}{\partial X^j} dX^j
\wedge dX^{i_1}\wedge\ldots\wedge dX^{i_k},
$$
\begin{multline*}
  [U,V]^{k_2\dots k_{g+\ell}} =
  \varepsilon^{k_2\dots k_{g+\ell}}_{i_2\dots i_g i_{g+1}\dots i_{g+\ell}}
U^{ri_2\dots i_g} \frac{\partial}{\partial X^r}
   V^{i_{g+1}\dots i_{g+\ell}}
%
%
 \\
+ \, (-1)^g
%
%
\varepsilon^{k_2\dots k_{g+\ell}}_{i_1\dots i_g i_{g+2}\dots i_{g+\ell}}
   V^{r i_{g+2}\dots i_{g+\ell}}
   \frac{\partial}{\partial X^r} U^{i_1\dots i_g}.
\end{multline*}

\subsection{Realizations of tensor operations}
\label{subsec-real}

   Let $\Lm$ be a finite dimensional free $\A$-module.
   Choose a basis $\{f_i\}$ in~$\Lm$ (over $\A$) and let $\{f^i\}$ be the dual basis in $\Lm^*$.
   For any $X\in \Lm$ we have $X=X^i f_i=x^{ia} f_i e_a$ where $x^{ia}\in \R$.
   Thus, the elements $f_{ia}:=f_ie_a$ form a basis of $\Lm$ considered as an $\R$-module.
   Let $f^{ia}=p\circ (f^i e^a):\Lm \to \R$.
   Then $\{f^{ia}\}$  form a basis of  $\Lm^*$ considered as an $\R$-module,
dual to  $\{f_{ia}\}$.
   In fact, $f^{ia}(f_{jb})=p(f^i e^a(f_j e_b))=p(e^a e_b) \delta^i_j =
\delta^{ia}_{jb}$ by (\ref{pee}) and $f^{ia}(X)=x^{ia}$.

   Let $t: \Lm\times \dots \times \Lm \to \A$ be an
$\A$-linear covariant tensor.

\begin{definition}
{\bf\cite{VSS}}
  The {\it realization} of $t$  is the $\R$-linear tensor
$$
R(t):=p\circ t: \Lm\times \dots \times \Lm \to \R.
$$
\end{definition}

   Let $t=t_{i_1\dots i_k} f^{i_1}\otimes_\A\ldots\otimes_\A f^{i_k}$.
   Then $t(X_1,\dots, X_k)=t_{i_1\dots i_k} X_1^{i_1} \dots X_k^{i_k}$, hence
$$
\begin{array}{l}
R(t)(X_1,\dots, X_k)=p(t_{i_1\dots i_k} X_1^{i_1} \dots X_k^{i_k})=\\[5pt]
\qquad\qquad
= p(t_{i_1\dots i_k} e_{a_1} \dots e_{a_k} x_1^{i_1a_1} \dots
x_k^{i_ka_k})= t_{i_1a_1\dots i_ka_k} x_1^{i_1a_1} \dots
x_k^{i_ka_k},
\end{array}
$$
where
\begin{equation}
\label{real-cov}
t_{i_1a_1\dots i_ka_k} :=
p(t_{i_1\dots i_k} e_{a_1}\dots e_{a_k}).
\end{equation}
  Thus
$$
R(t)=t_{i_1a_1\dots i_ka_k}
f^{i_1a_1}\otimes_\R\ldots \otimes_\R f^{i_ka_k},
$$
and its components $t_{i_1a_1\dots i_ka_k}$ can be calculated
by~(\ref{real-cov}).

   Let $\Lm'$ be another finite-dimensional $\A$-module and $\Psi: \Lm'\to \Lm$ be an
$\A$-linear map.
   Let us show that
\begin{equation}
\label{Alin-map-real}
  \Psi^* R(t) = R(\Psi^* t).
\end{equation}

    In fact, let $\{g_\alpha\}$ be a basis in $\Lm'$.
    Denote the components of $\Psi$ by
$\Psi_\alpha^i=\Psi_\alpha^{ic} e_c$.
    Then the components $\Psi_{\alpha b}^{ia}$ of $\Psi$ considered as an
$\R$-linear map $\Lm'\to \Lm$ are
$\Psi_{\alpha b}^{ia} = \Psi_\alpha^{ic} \gamma_{bc}^a$.
    We have
$$
(\Psi^* t)_{\alpha_1\dots \alpha_k} =
\Psi_{\alpha_1}^{i_1} \dots \Psi_{\alpha_k}^{i_k} t_{i_1\dots i_k},
$$
whence
\begin{align*}
  (R(\Psi^* t))_{\alpha_1b_1\dots \alpha_kb_k} &
 = p(\Psi_{\alpha_1}^{i_1} \dots \Psi_{\alpha_k}^{i_k}
t_{i_1\dots i_k} e_{b_1}\dots e_{b_k}) \\[10pt]
 &= p(\Psi_{\alpha_1}^{i_1c_1} e_{c_1}e_{b_1}\dots
\Psi_{\alpha_k}^{i_kc_k} e_{c_k} e_{b_k}
t_{i_1\dots i_k})                          \\[10pt]
 &= \Psi_{\alpha_1}^{i_1 c_1} \gamma_{c_1b_1}^{a_1}\dots
\Psi_{\alpha_k}^{i_k c_k} \gamma_{c_kb_k}^{a_k}
p(t_{i_1\dots i_k} e_{a_1}\dots e_{a_k})  \\[10pt]
 &= \Psi_{\alpha_1b_1}^{i_1a_1}\dots
\Psi_{\alpha_kb_k}^{i_ka_k}
R(t)_{i_1a_1\dots i_ka_k} =
(\Psi^* R(t))_{\alpha_1b_1\dots \alpha_kb_k}.\\
\end{align*}


  The realization of 1-forms $R:\omega\mapsto R(\omega)$ is an $\R$-linear
isomorphism from $\Lm^*$ considered as a vector space onto the dual vector space
$\Lm_\R^*$ to $\Lm$ (considered as a vector space over $\R$).
  The isomorphism $R$ transfers the structure of an $\A$-module from $\Lm^*$ to
$\Lm_\R^*$: if $\xi=R(\omega)$, $\alpha\in \A$, then $\alpha \xi= R(\alpha \omega)$.
  The structure of an $\A$-module on $\Lm_\R^*$ can be described as follows: if $v\in\L$,
then $(\alpha \xi) (v) = \xi(\alpha v)$.
  In fact, $(\alpha \xi) (v)= R(\alpha \omega)(v) = p\circ \alpha \omega(v) =
p\circ \omega(\alpha v) = R(\omega)(\alpha v) = \xi(\alpha v)$.
  It will be convenient in the sequel to identify the modules $\Lm^*$ and $\Lm_\R^*$ and
consider $\Lm^*$ as a dual $\A$-module to $\Lm$ with  contraction
$\Lm\times \Lm^* \ni (\omega, v) \mapsto \langle \omega, v\rangle_\A \in \A$ and as
a dual vector space to $\Lm$ with contraction
$ \Lm\times \Lm^* \ni (\omega, v) \mapsto \langle \omega, v\rangle_\R = p\circ\langle \omega, v\rangle_\A \in \R$.


   Now we describe the realization of contravariant tensors.

   Let $\Lm$ be a finite-dimensional $\A$-module and
$u: \Lm^*\times \dots \times \Lm^* \to \A$ be
an $\A$-linear contravariant tensor.
\begin{definition}
   The {\it realization} of  $u$ is the $\R$-linear tensor
$$
R(u):=p\circ u: \Lm^*\times \dots \times \Lm^* \to \R.
$$
\end{definition}

  Making use of the diagram
$$
\xymatrix{%
{\Lm^*\times \ldots \times \Lm^*} \ar[r]^-{u} &
{\A} \ar[rr]^{\varphi} \ar[dr]_{p} & & {\A^*} \ar[dl]^{1_\A}\\
 & & {\R} \\
}
$$
we may represent the realization $R(u)$ as  $1_\A \circ \varphi\circ u$, where
$\varphi:\A\to \A^*$ is  isomorphism~(\ref{phi}) induced by the
Frobenius form~$q$.

   Let $u=u^{i_1\dots i_k} f_{i_1}\otimes_\A\ldots\otimes_\A f_{i_k}$ be the representation of $u$
in terms of a basis $\{f_i\}$ in $\Lm$, and let
$u^{i_1\dots i_k} e^{a_1}\dots e^{a_k} =
u^{i_1a_1\dots i_ka_k}_b e^b$ be the expansion in terms of the basis of $\A$.
   From (\ref{F-contr}) we have
$p(u^{i_1\dots i_k} e^{a_1}\dots e^{a_k})=
u^{i_1a_1\dots i_ka_k}_b \delta^b$.
   Denote
\begin{equation}
\label{real-contrav}
u^{i_1a_1\dots i_ka_k} :=
p(u^{i_1\dots i_k} e^{a_1}\dots e^{a_k}).
\end{equation}
   Then
$$
R(u)=u^{i_1a_1\dots i_ka_k}
f_{i_1a_1}\otimes_\R\ldots \otimes_\R f_{i_ka_k}.
$$

   Let now $\Lm'$ be another finite-dimensional $\A$-module,
and let $\Psi: \Lm\to \Lm'$ be an $\A$-linear map.
   Suppose that  a contravariant $\A$-linear
tensor $u$ on $\Lm$  is   $\Psi$-related to a contravariant $\A$-linear tensor
$v$ on $\Lm'$.

\begin{proposition}
\label{prop-real-contrav}
   $R(u)$ is  $\Psi$-related to  $R(v)$.
\end{proposition}

\begin{Proof}
   The condition that $u$ is $\Psi$-related to $v$ is as follows~\cite{KMS, daS-W}
$$
v^{\alpha_1\dots \alpha_k} =
\Psi^{\alpha_1}_{i_1} \dots \Psi^{\alpha_k}_{i_k} u^{i_1\dots i_k}.
$$
   We have
\begin{multline*}
  (R(v))^{\alpha_1b_1\dots \alpha_kb_k}
= p(\Psi^{\alpha_1}_{i_1} \dots \Psi^{\alpha_k}_{i_k} u^{i_1\dots
i_k} e^{b_1}\dots e^{b_k}) = p(\Psi^{\alpha_1c_1}_{i_1}
e_{c_1}e^{b_1}\dots
\Psi^{\alpha_kc_k}_{i_k} e_{c_k} e^{b_k}
u^{i_1\dots i_k}).
\end{multline*}

  Note that $e^{b} e_{c} = \gamma_{ac}^b e^a$ by~(\ref{ee-gamma}).
  Thus
\begin{multline*}
  (R(v))^{\alpha_1b_1\dots \alpha_kb_k}
 = \Psi^{\alpha_1c_1}_{i_1} \gamma_{a_1c_1}^{b_1}\dots
\Psi^{\alpha_kc_k}_{i_k} \gamma_{a_kc_k}^{b_k}
p(u^{i_1\dots i_k} e^{a_1}\dots e^{a_k}) \\ =
\Psi^{\alpha_1b_1}_{i_1a_1}\dots
\Psi^{\alpha_kb_k}_{i_ka_k}
(R(u))^{i_1a_1\dots i_ka_k}.
\end{multline*}
\end{Proof}

  Let us find the  expression for $R(v)$ where $v\in \Lm$.
  For a basis $\{f_i\}$  in  $\Lm$ we have $v=v^j f_j$.
  Let $v^j=v^{jb} e_b = v^j_b e^b$.
  Then $v^j e^a = v^{jb} e_b e^a = v^j_b e^a e^b$.
  Since $p(e_a e^b) =\delta_a^b$ and $p(e^a e^b)=q^{ab}$
by  $(\ref{pee})$ and $(\ref{pee-up})$, respectively, we get
\begin{equation}
\label{real-vector}
(R(v))^{ja}=v^{ja} = v^j_b q^{ab}.
\end{equation}

\begin{proposition}
\label{prop-real-cont}
  Let  $\Lm$ be a finite-dimensional $\A$-module, $v\in \Lm$, and let
$t\in {T}^{k,0}(\Lm)$ be a covariant $\A$-tensor.
  Then
\begin{equation}
\label{real-cont}
   R(i(v) t) = i(R(v)) R(t),
\end{equation}
where $i(v)$, $i(R(v))$ are defined by~{\rm (\ref{i(u)})}.
\end{proposition}

\begin{Proof}
  Let $\{f_i\}$ be a basis in  $\Lm$,
$v=v^j f_j$, and
$t=t_{i_1\dots i_k} f^{i_1} \otimes \ldots \otimes f^{i_k}$.
  Denote $\theta=i(v)t$.
  We compute the components of $R(v)$, $R(t)$ and $R(\theta)$.

  We have
$(R(v))^{ja} = v^j_b q^{ab}$
by $(\ref{real-vector})$.

  We define $\gamma_{a_1\dots a_k}^b\in\R$ by
\begin{equation}
\label{gammaup}
e_{a_1}\dots e_{a_k} = \gamma_{a_1\dots a_k}^b e_b.
\end{equation}
   Clearly, $\gamma_{a_1\dots a_k}^b =
\gamma_{a_1a_2}^{c_1}\gamma_{c_1a_3}^{c_2} \dots
\gamma_{c_{k-1}a_k}^b$ (see~\cite{Kruch, VSS}).

  Let $t_{ji_1\dots i_{k-1}} = t_{ji_1\dots i_{k-1}}^s e_s$ be the expansion in terms of
the basis in $\A$.
  Then
$$
t_{ji_1\dots i_{k-1}} e_a e_{a_1}\dots e_{a_{k-1}} = t_{ji_1\dots
i_{k-1}}^s e_s e_a e_{a_1}\dots e_{a_{k-1}} = t_{ji_1\dots
i_{k-1}}^s
\gamma_{sa{a_1}\dots a_{k-1}}^c e_c.
$$
  Contracting with $p$ we obtain
$$
(R(t))_{jai_1a_1\dots i_{k-1}a_{k-1}}=
t_{jai_1a_1\dots i_{k-1}a_{k-1}}=
t_{ji_1\dots i_{k-1}}^s
\gamma_{sa{a_1}\dots a_{k-1}}^c p_c.
$$

  We also have
$$
\theta_{i_1\dots i_{k-1}} =
v^j t_{ji_1\dots i_{k-1}} = v^i_b e^b t_{ji_1\dots i_{k-1}}^s e_s = v^i_b
q^{bd} t_{ji_1\dots i_{k-1}}^s e_d e_s = v^i_b q^{bd} t_{ji_1\dots i_{k-1}}^s
\gamma_{sd}^r  e_r
$$
  and
\begin{multline*}
\theta_{i_1\dots i_{k-1}} e_{a_1}\dots e_{a_{k-1}} =
v^i_b t_{ji_1\dots i_{k-1}}^s q^{bd} \gamma_{sd}^r
e_r e_{a_1}\dots e_{a_{k-1}} \\[5pt]
 =
v^i_b t_{ji_1\dots i_{k-1}}^s q^{bd} \gamma_{sd}^r
\gamma_{r {a_1}\dots a_{k-1}}^c e_c=
v^i_b t_{ji_1\dots i_{k-1}}^s q^{bd}
\gamma_{sd{a_1}\dots a_{k-1}}^c e_c.
\end{multline*}
  Hence
$$
(R(\theta))_{i_1a_1\dots i_{k-1}a_{k-1}}=
v^i_b t_{ji_1\dots i_{k-1}}^s q^{bd}
\gamma_{sd {a_1}\dots a_{k-1}}^c p_c,
$$
which means that
\begin{multline*}
(i(R(v)) R(t))_{i_1a_1\dots i_{k-1}a_{k-1}} = v^{ja}
t_{jai_1a_1\dots i_{k-1}a_{k-1}}\\[5pt] =
 v^j_b q^{ab} t_{ji_1\dots i_{k-1}}^s
 \gamma_{sa{a_1}\dots a_{k-1}}^c p_c =
  (R(\theta))_{i_1a_1\dots i_{k-1}a_{k-1}}.
\end{multline*}
\end{Proof}

\begin{remark}
\label{real-inj}
   Note that for a Frobenius Weil algebra $\A$ the realization of
tensors is an injective operation.
   Indeed, if $t:\Lm\times \ldots \times \Lm \to \A$ is a covariant
$\A$-tensor then
$R(t) = p\circ t = q(t,1_\A)$.
  Since $q$ is nondegenerate, it follows that
$t$ and $R(t)$ vanish or do not vanish simultaneously.
  For the contravariant $\A$-tensors proof is similar.
\end{remark}

\begin{remark}
\label{rem-mixed}
  The realization of a tensor $t$ of type $(k,\ell)$ for $k,\ell\ge1$
can be constructed in the same way.
  Namely, if
$$
t=t_{i_1\dots i_k}^{j_1\dots j_\ell}
f^{i_1}\otimes_\A\ldots\otimes_\A f^{i_k} \otimes_\A
f_{j_1}\otimes_\A\ldots\otimes_\A f_{j_\ell}
$$
then
$$
R(u)=t_{i_1a_1\dots i_ka_k}^{j_1b_1\dots j_\ell b_\ell}
f^{i_1a_1}\otimes_\R\ldots \otimes_\R f^{i_ka_k}\otimes_\R
f_{j_1b_1}\otimes_\R\ldots \otimes_\R f_{j_\ell b_\ell}.
$$
where
$$
t_{i_1a_1\dots i_ka_k}^{j_1b_1\dots j_\ell b_\ell} :=
p(t_{i_1\dots i_k}^{j_1\dots j_\ell} e_{a_1}\dots e_{a_k} e^{b_1}\dots e^{b_\ell}) .
$$
\end{remark}


\subsection{The complete lift of a covariant tensor field}
\label{subsec-cl-cov}

  Let $M$ be a smooth manifold of dimension~$m$ and $\pi_\A:T^\A M\to M$
 its Weil bundle.
  For a local chart $(U,x^1,\dots, x^m)$ on $M$ the functions
$X^i=(x^i)^\A=x^{ia}e_a$ form a system of $\A$-valued local coordinates on
$T^\A U\subset T^\A M$, and $(x^{ia})$ are real
local coordinates on~$T^\A U$, $x^{i0}=x^i\circ \pi_\A$.

  Let  $\xi\in {T}^{k,0}(M)$ be a tensor field of type~$(k,0)$ on~$M$.
  In local coordinates
$$
\xi=\xi_{i_1\dots i_k}
dx^{i_1}\otimes\ldots\otimes dx^{i_k}.
$$
  Consider the analytic prolongations
$\Xi_{i_1\dots i_k}=(\xi_{i_1\dots i_k})^\A$ of the functions
$\xi_{i_1\dots i_k}$.
  The analytic prolongation
$\xi^\A\in {T}^{k,0}_{\Adiff}(T^\A M)$ of~$\xi$
locally is of the form
$\xi^\A= \Xi_{i_1\dots i_k} dX^{i_1}\otimes\ldots\otimes dX^{i_k}$.
  Denote
\begin{equation}
\label{lift-cov}
\xi_{i_1a_1\dots i_ka_k} :=
p(\Xi_{i_1\dots i_k} e_{a_1}\dots e_{a_k}).
\end{equation}
   We have
$$
R(\xi^\A)= \xi_{i_1a_1\dots i_ka_k}
dx^{i_1a_1}\otimes\ldots\otimes dx^{i_ka_k}.
$$

\begin{definition}
{\bf\cite{ELOS, Shir12}}
  The {\it complete lift} of~$\xi$
is the tensor field
$$
\xi^C=\xi^C_\A:=R(\xi^\A)
$$
on $T^\A M$.
\end{definition}

   From (\ref{Alin-map-real}) it follows that
for every smooth map~$\varphi : N\to M$
$$
(T^\A \varphi)^* (\xi^C)  =(\varphi^* \xi)^C.
$$

\begin{remark}
\label{cl-inj-cov}
  It follows from the Remark~\ref{real-inj} that for $k\ge 1$
the complete lift induces the injective map
${T}^{k,0}(M)\to {T}^{k,0}(T^\A M)$, that is,
$\xi^C = 0$ if and only if $\xi=0$.
\end{remark}

\begin{proposition}
\label{prop-inj-cov}
   The complete lift is an injective map
$C^\infty(M)\to C^\infty(T^\A M)$ if and only if
$p(1_\A)=p_0\ne0$. When $p(1_\A)=0$ its kernel is the space of locally constant functions.
\end{proposition}

\begin{Proof}
   In fact, let $f\in C^\infty(M)$ be a non-zero function
and $f^\A = f^a e_a=f^0e_0+f^{\hat a}e_{\hat a}$. Then $f^C=p_a f^a$ (recall
that the indices $\hat a,
\hat b$ run through the set of values 1, \dots,
$n$).
   It follows from~(\ref{A-prol}) that
the functions $f^{\hat a}$ are locally of the form
$$
  f^{\hat a} = \frac{\partial f}{\partial x^i} x^{ia} \, + \,
    \hbox{summands of degree $\ge2$ in $x^{j\hat b}$ }.
$$

   Thus $f^{\hat a}\equiv 0$ if and only if all partial derivatives
$\frac{\partial f}{\partial x^i}$ vanish, that is, if
$f\equiv{\rm const}$.
   This means that the condition $f^C=p_0f^0 + p_{\hat a}f^{\hat a} = 0$
is equivalent to $p_0f^0=p_{\hat a}f^{\hat a}=0$.
   When  $p(1_\A)\ne 0$ the first equation gives $f=0$.
   In case $p(1_\A)=0$ the second equation, by virtue of the fact that $p_n\ne0$,
means that $f^n=0$. Consequently, $f={\rm const}$.
\end{Proof}

\begin{proposition}
\label{prop-ext-dif}
  Let $M_\A$ be an $\A$-smooth manifold and let
$\Xi\in\Omega^k_{\Adiff} (M_\A)$ be an $\A$-smooth exterior form.
  Then
\begin{equation}
\label{real-ext-dif}
  R(d\,\Xi)=   d\, (R(\Xi)),
\end{equation}
where $d$ is the exterior differential.
\end{proposition}

\begin{Proof}
  If
$$
\Xi=\Xi_{i_1\dots i_k}
dX^{i_1}\wedge\ldots\wedge dX^{i_k}
$$
in terms of local coordinates, then
$$
d\,\Xi=\frac{\partial \Xi_{i_1\dots i_k}}{\partial X^j}
dX^j \wedge dX^{i_1}\wedge\ldots\wedge dX^{i_k}.
$$

  Let  $\Xi_{i_1\dots i_k} e_{a_1} \dots e_{a_k} =
\Xi_{i_1a_1\dots i_k a_k}^c e_c$ be expansions in terms of the basis of $\A$.
  Then
$$
(R(\Xi))_{i_1a_1\dots i_ka_k}=\Xi_{i_1a_1\dots i_k a_k}^c p_c.
$$
  From~(\ref{u-s3}) we obtain
\begin{multline}
\label{dRXi}
(d\,R(\Xi))_{jbi_1a_1\dots i_ka_k}=
\frac{\partial}{\partial x^{jb}}
  (R(\Xi))_{i_1a_1\dots i_ka_k}=
\frac{\partial}{\partial x^{jb}}
  \Xi_{i_1a_1\dots i_k a_k}^c\,p_c  \\[10pt]
=
\gamma^c_{bd} \delta^g \frac{\partial}{\partial x^{jg}}
\Xi_{i_1a_1\dots i_k a_k}^d\,p_c =
q_{bd} \delta^g \frac{\partial}{\partial x^{jg}}
\Xi_{i_1a_1\dots i_k a_k}^d.
\end{multline}

  On the other hand, by (\ref{der-dir}),
$$
({d\,\Xi})_{ji_1\dots i_k} e_{b} e_{a_1}\dots e_{a_k} =
\frac{\partial}{\partial X^j}
\Xi_{i_1a_1\dots i_k a_k}^c\, e_c e_b=
\delta^g \frac{\partial}{\partial x^{jg}}
    \Xi_{i_1a_1\dots i_k a_k}^c\, \gamma^d_{bc} e_d.
$$
  Therefore
\begin{equation}
\label{RdXi}
R(d\Xi)_{jbi_1a_1\dots i_ka_k}=
\delta^g \frac{\partial}{\partial x^{jg}}
    \Xi_{i_1a_1\dots i_k a_k}^c\,
    \gamma^d_{bc} p_d =
\delta^g \frac{\partial}{\partial x^{jg}}
    \Xi_{i_1a_1\dots i_k a_k}^c\, q_{bc},
\end{equation}
which coincides with~(\ref{dRXi}).
\end{Proof}

\begin{corollary}
\label{cor-ext-dif}
  Let $M$ be a smooth manifold and $T^\A M$ be its Weil bundle.
  The complete lift commutes with the exterior differential, i.e.
\begin{equation}
\label{lift-ext-dif}
    (d\,\xi)^C = d(\xi^C)
\end{equation}
for every $\xi \in \Omega^*(M)$.
\end{corollary}

   This means that the complete lift of exterior forms induces
a homomorphism of de~Rham cohomology spaces
\begin{equation}
\label{lift-dR}
     H^*_{dR}(M) \longrightarrow H^*_{dR}(T^\A M), \qquad
        [\xi] \longmapsto [\xi^C].
\end{equation}

\begin{theorem}
\label{cl-dR}
  Let $(\A,q)$ be a Frobenius Weil algebra and $p$  its Frobenius covector.
  If $p(1_\A)\ne 0$, then the homomorphism~{\rm(\ref{lift-dR})} is an isomorphism.
  If $p(1_\A)=0$, then
the homomorphism~{\rm(\ref{lift-dR})} is zero.
\end{theorem}

\begin{Proof}
    The manifold $M$ may be embedded into $T^\A M$ by means of the zero section
$$
s_\A: M\to T^\A M.
$$
   It is shown in~\cite{VSh96} that  the complexes
$(\Omega^*_{\Adiff} (T^\A M), d)$ and $\A  \otimes (\Omega^*(M),d)$ are isomorphic.
   In fact, to any  $\A$-valued exterior form $\xi$ on~$M$
there corresponds its $\A$-prolon\-ga\-tion $\xi^{\A}$
which is an $\A$-smooth form on $T^{\A}M$.
   Moreover, each $\A$-smooth form $\theta$ on $T^\A M$ coincides with
the $\A$-prolon\-ga\-tion of $\theta|_{s_\A(M)}=s_\A^*(\theta)$.
   Denote the cohomology of $(\Omega^*_{\Adiff} (T^\A M),d)$ by
$H^*_{\Adiff} (T^\A M)$.
   Thus,
$$
H^*_{\Adiff} (T^\A M) \cong \A  \otimes H^*_{dR}(M).
$$

   It is also clear that the maps
$\pi_\A^* : H^*_{dR}(M)\to H^*_{dR}(T^\A M)$ and
$s_\A^*: H^*_{dR}(T^\A M) \to H^*_{dR}(M)$ are mutually inverse isomorphisms
of de~Rham cohomology (we use here the symbols
$\pi^*_\A$ and $s^*_\A$ simultaneously for the maps  of exterior forms and for the corresponding maps of the de~Rham cohomology).
   Therefore we get the isomorphism
$$
  H^*_{\Adiff} (T^\A M) \;
  \stackrel{s_\A^*}{\hbox to 12mm{\rightarrowfill}} \;
  \A  \otimes H^*_{dR}(M) \;
  \stackrel{\A\otimes \pi_\A^*}{\hbox to 12mm{\rightarrowfill}}
   \; \A \otimes H^*_{dR}(T^\A M).
$$

   Let $\xi\in \Omega^*(M) = \R\otimes \Omega^*(M)
\subset \A\otimes \Omega^*(M)$ be a closed form and
$\xi^\A = \xi^0 e_0 +\ldots + \xi^n e_n \in
\Omega^*_{\Adiff}(T^\A M)$ its analytic prolongation.
   It  is easily seen that
$$
s^*_\A: \Omega^*_{\Adiff}(T^\A M) \to \A \otimes \Omega^*(M)
$$
maps $\xi^0$ to $\xi$ and all the forms  $\xi^1$, \dots, $\xi^n$ to zero.
  Consequently,
$$
\pi_\A^* \circ s^*_\A: H^*_{\Adiff} (T^\A M) \to
\A \otimes H^*_{dR}(T^\A M)
$$
maps the cohomology class $[\xi^0]$ to $1_\A\otimes[\xi^0]$ and
the classes $[\xi^1], \dots, [\xi^n]$ to zero.
   It follows from~(\ref{lift-cov}) that
$\xi^C = p_a \xi^a = p_0\xi^0 + p_1\xi^1+
\ldots + p_n \xi^n$.
   Hence $[\xi^C] = p_0 [\xi^0] = p_0[\pi^*_\A\xi]\in
H^*_{dR}(T^\A M)$.
\end{Proof}

  Thus, each Frobenius Weil algebra $(\A,q)$ determines the endomorphism of the
coho\-mo\-logy spaces
$$
  \Delta_{\A,q} : H^*_{dR}(M)\;
\overset{C}{\hbox to 12mm{\rightarrowfill}}
\; H^*_{dR}(T^\A M) \;
\overset{(\pi_\A^*)^{-1}}{\hbox to 12mm{\rightarrowfill}}
\; H^*_{dR}(M),
$$
where the first arrow means the complete lift.
  It follows from the Theorem~\ref{cl-dR} that this endomorphism is
the multiplication by  $p(1_\A)$.

\begin{example}
\label{ex-tang-b}
  Let
$\tau = \pi_{\R(\varepsilon)} : TM\to M$ be the tangent bundle.
  We denote the local coordinates on~$M$ by $(x^i)$ and the
corresponding local coordinates on~$TM$ by~$(x^i, y^i)$.

  Take the basis  $\{e_0=1, e_1=\varepsilon\}$ in~$\R(\varepsilon)$, and
introduce two Frobenius covectors $p_{(0)}$ and
$p_{(1)}$ on $\R(\varepsilon)$ defined by $p_{(0)}(1)=0$, $p_{(0)}(\varepsilon)=1$ and
$p_{(1)}(1)=1$, $p_{(1)}(\varepsilon)=1$, respectively.

  Let
$\xi=\xi_{i_1\dots i_k} dx^{i_1}\wedge \ldots
\wedge dx^{i_k}\in \Omega^k(M)$ be a closed form.
  For  $p_{(0)}$ the corresponding complete lift $\xi^C_{(0)}$ is
$$
  \xi^C_{(0)} =
  y^j \frac{\partial \xi_{i_1\dots i_k}}{\partial x^j}
   dx^{i_1}\wedge \ldots \wedge dx^{i_k} +
   \xi_{i_1\dots i_k} dx^{i_1}\wedge \ldots
    \wedge dx^{i_{k-1}} \wedge dy^{i_k}.
$$
  It is easily seen that
$$
\eta:=y^j \xi_{j i_1\dots i_{k-1}}
dx^{i_1}\wedge\ldots\wedge dx^{i_{k-1}}
$$
is a well-defined form on~$TM$ and that
$\xi$ is closed if and only if $\xi^C_{(0)} = d\eta$.
Thus,
$[\xi^C_{(0)}]=0\in H^*_{dR}(TM)$.

  For $p_{(1)}$ the complete lift $\xi^C_{(1)}$ is
$$
\begin{array}{l}
  \xi^C_{(1)} =
   \xi_{i_1\dots i_k}  dx^{i_1}\wedge \ldots \wedge dx^{i_k}~ + \\[5pt]
  \qquad +~ y^j \frac{\partial \xi_{i_1\dots i_k}}{\partial x^j}
   dx^{i_1}\wedge \ldots \wedge dx^{i_k} +
   \xi_{i_1\dots i_k} dx^{i_1}\wedge \ldots
    \wedge dx^{i_{k-1}} \wedge dy^{i_k} = \\[5pt]
    \qquad = \tau^*\xi + \xi^C_{(0)}.
\end{array}
$$
  Therefore $[\xi^C_{(1)}]= [\tau^*\xi]\in H^*_{dR}(TM)$.
\end{example}


\begin{remark}
   Let  $\xi\in {T}^{k,0}(M)$ be a tensor field of type
$(k,0)$ on~$M$ and
$\xi^\A \in {T}^{k,0}_{\Adiff}(T^\A M)$  its
analytic prolongation to $T^\A M$.
   Then $e_n \xi^\A$ is also  an $\A$-smooth tensor field.
   If
$\xi=\xi_{i_1\dots i_k}
dx^{i_1}\otimes\ldots\otimes dx^{i_k}$, then
$e_n \xi^\A= e_n \xi_{i_1\dots i_k} dX^{i_1}\otimes\ldots\otimes dX^{i_k}
= e_n \xi_{i_1\dots i_k} dx^{i_10}\otimes\ldots\otimes dx^{i_k0}$.
   It follows that
$R(e_n \xi^\A)\in {T}^{k,0}(T^\A M)$ is of the form
$$
 R(e_n \xi^\A)= \xi_{i_1\dots i_k} dx^{i_10}\otimes\ldots\otimes dx^{i_k0},
$$
that is, it coincides with $\pi_\A^* \xi$.
   Thus,
\begin{equation}
\label{xiV}
  \pi_\A^* \xi = R(e_n \xi^\A).
\end{equation}
\end{remark}


\subsection{The complete lift of a contravariant tensor field}
\label{subsec-cl-contrav}

  Let $M$  be a real smooth manifold,
$T^\A M$  its Weil bundle, and let
$u\in{T}^{0,k}(M)$ be a contravariant tensor field on~$M$.
  In terms of local coordinates,
$$
u=u^{i_1\dots i_k} \frac{\partial}{\partial x^{i_1}}
\otimes\ldots\otimes \frac{\partial}{\partial x^{i_k}}.
$$
  Consider the analytic prolongations~$U^{i_1\dots i_k}=(u^{i_1\dots i_k})^\A$
of the functions $u^{i_1\dots i_k}$.
  The analytic prolongation of~$u$ is
$$
u^\A=U^{i_1\dots i_k} \frac{\partial}{\partial X^{i_1}}
\otimes\ldots\otimes \frac{\partial}{\partial X^{i_k}}.
$$
  Let $U^{i_1\dots i_k} e^{a_1}\dots e^{a_k} =
U^{i_1a_1\dots i_ka_k}_b e^b$ be the expansions in terms of the basis in $\A$.

  Denote
\begin{equation}
\label{lift-contrav}
u^{i_1a_1\dots i_ka_k} :=
p(U^{i_1\dots i_k} e^{a_1}\dots e^{a_k}) =
U^{i_1a_1\dots i_ka_k}_b \delta^b.
\end{equation}
  We have
$$
   R(u^\A) = u^{i_1a_1\dots i_ka_k}
\frac{\partial}{\partial x^{i_1a_1}}\otimes\ldots\otimes
\frac{\partial}{\partial x^{i_ka_k}}.
$$

\begin{definition}
  The {\it complete lift} $u^C$ of $u$ is the tensor field
$$
u^C=u^C_\A := R(u^\A)
$$
on $T^\A M$.
\end{definition}

   It follows immediately from the Proposition~\ref{prop-real-contrav} that
if $\varphi : M\to N$ is  a smooth map and a tensor field $u$ is
$\varphi$-related with a tensor field $v$ on $N$, then
$u^C$ is $T^\A \varphi$-related to~$v^C$.

\begin{remark}
\label{cl-inj}
  It follows from the Remark~\ref{real-inj} that
for $k\ge 1$ the complete lift is an injective map
${T}^{0,k}(M)\to {T}^{0,k}(T^\A M)$.
\end{remark}


\begin{proposition}
\label{prop-Sch-br}
  Let $M_\A$ be an  $\A$-smooth manifold and let
$U, V \in {\cal V}^*_{\Adiff}(M_\A)$ be two
$\A$-smooth multivector fields.
  Then
\begin{equation}
\label{real-Sch-br}
  [R(U), R(V)] = R([U,V]).
\end{equation}
\end{proposition}

\begin{Proof}
  Let $U\in {\cal V}_{\Adiff}^g(M_\A)$,
$V\in{\cal V}^\ell_{\Adiff}(M_\A)$.
  According to~(\ref{Sch-br}), in terms of local coordinates, the Schouten-Nijenhuis bracket
$[U,V]$ is of the form
\begin{multline}
\label{Sch-br-an}
  [U,V]^{k_2\dots k_{g+\ell}} =
%
%
\varepsilon^{k_2\dots k_{g+\ell}}_{i_2\dots i_g i_{g+1}\dots i_{g+\ell}}
U^{ri_2\dots i_g} \frac{\partial}{\partial X^r}
   V^{i_{g+1}\dots i_{g+\ell}}  \\[5pt]
+ \, (-1)^g
%
%
\varepsilon^{k_2\dots k_{g+\ell}}_{i_2\dots i_{g+1} i_{g+2}\dots i_{g+\ell}}
   V^{r i_{g+2}\dots i_{g+\ell}}
   \frac{\partial}{\partial X^r} U^{i_2\dots i_{g+1}}.
\end{multline}

  Let us multiply each side of (\ref{Sch-br-an}) by
$e^{a_2}\dots e^{a_{g+\ell}}$ and then contract with $\delta^s$.
  In the left-hand side we get
$(R([U,V]))^{k_2a_2\dots k_{g+\ell}a_{g+\ell}}$.
  Using (\ref{der-dir}) and (\ref{partialF}), we transform  the first summand in the right-hand side of (\ref{Sch-br-an}), omitting
the coefficient $\varepsilon^{k_2\dots k_{g+\ell}}_{i_2\dots i_g i_{g+1}\dots
i_{g+\ell}}$ as follows.
\begin{multline*}
U^{ri_2\dots i_g} e^{a_2}\dots e^{a_g}
\frac{\partial}{\partial X^r} V^{i_{g+1}\dots i_{g+\ell}}
e^{a_{g+1}}\dots e^{a_{g+\ell}}\\[10pt]
 =
U^{ri_2\dots i_g} e^{a_2}\dots e^{a_g}
\delta^b\frac{\partial}{\partial x^{rb}}
V^{i_{g+1}a_{g+1}\dots i_{g+\ell}a_{g+\ell}}
\\[10pt]
 =
U^{ri_2\dots i_g} e^a e^{a_2}\dots e^{a_g}
\delta^b\frac{\partial}{\partial x^{rb}}
V^{i_{g+1}a_{g+1}\dots i_{g+\ell}a_{g+\ell}}_a \\[10pt]
 =
U^{ri_2\dots i_g} e^a e^{a_2}\dots e^{a_g}
\frac{\partial}{\partial x^{ra}}
(\delta^b V^{i_{g+1}a_{g+1}\dots i_{g+\ell}a_{g+\ell}}_b)
 \\[10pt]
=
U_s^{ra\,i_2a_2\dots i_ga_g} e^s
\frac{\partial}{\partial x^{ra}}
((R(V))^{i_{g+1}a_{g+1}\dots i_{g+\ell}a_{g+\ell}}).
\end{multline*}
  Contracting the result with $\delta^s$ we get
$$
\displaystyle
(R(U))^{ra\,i_2a_2\dots i_ga_g}
\frac{\partial}{\partial x^{ra}}
\Bigl((R(V))^{i_{g+1}a_{g+1}\dots i_{g+\ell}a_{g+\ell}}\Bigr).
$$
  The commutativity of multiplication in $\A$ yields that
$\varepsilon^{k_2\dots k_{g+\ell}}_{i_2\dots i_{g+\ell}}=
\varepsilon^{k_2a_2\dots k_{g+\ell}a_{g+\ell}}_
{i_2a_2\dots  i_{g+\ell}a_{g+\ell}}$.
  In the same manner we transform the second summand in
(\ref{Sch-br-an}).
  As a result, we have
\begin{multline*}
(R([U,V]))^{k_2a_2\dots k_{g+\ell}a_{g+\ell}} \\[15pt]
\displaystyle
=
\varepsilon^{k_2a_2\dots k_{g+\ell}a_{g+\ell}}_
{i_2a_2\dots i_ga_g i_{g+1}a_{g+1}\dots i_{g+\ell}a_{g+\ell}}
(R(U))^{ra\,i_2a_2\dots i_ga_g}
\textstyle\frac{\partial}{\partial x^{ra}}
\Bigl((R(V))^{i_{g+1}a_{g+1}\dots i_{g+\ell}a_{g+\ell}}\Bigr) \\[15pt]
 + \, (-1)^g
\varepsilon^{k_2a_2\dots k_{g+\ell}a_{g+\ell}}_
{i_1a_1\dots i_ga_g i_{g+2}a_{g+2}\dots i_{g+\ell}a_{g+\ell}}
(R(V))^{ra\, i_{g+2}a_{g+2}\dots i_{g+\ell}a_{g+\ell}}
\textstyle\frac{\partial}{\partial x^{ra}}
\Bigl((R(U))^{i_1a_1\dots i_ga_g}\Bigr),
\end{multline*}
which coincides with (\ref{real-Sch-br}).
\end{Proof}

\begin{corollary}
  Let $M$ be a smooth manifold and  $T^\A M$  its Weil bundle.
  The complete lift commutes with the Schouten-Nijenhuis
bracket, i.e.,
\begin{equation}
\label{lift-Sch-br}
  [u,v]^C = [u^C, v^C].
\end{equation}
for every $u,v \in {\cal V}^*(M)$.
\end{corollary}

\begin{proposition}
  Let $M$ be a smooth manifold, and let $\xi\in \Omega^*(M)$,
$v\in{\cal V}^1(M)$.
  Then
\begin{equation}
\label{lift-cont}
  i(v^C) \xi^C = (i(v)\xi)^C.
\end{equation}
\end{proposition}

\begin{Proof}
   Follows from the Proposition~\ref{prop-real-cont}.
\end{Proof}

\begin{remark}
  For $v\in{\cal V}^\ell(M)$, $\ell\ge 2$,
relation~(\ref{lift-cont}), in general, does not remain valid.
  Consider,  for example,  $\xi=\xi_{ij}dx^i\wedge dx^j\in \Omega^2(M)$,
$v=v^{ij}  \frac{\partial}{\partial x^i}\wedge
\frac{\partial}{\partial x^j}\in{\cal V}^2(M)$.
  Their complete lifts to the tangent bundle $TM$ corresponding to the
Frobenius covector~$p_{(0)}$ (see Example \ref{ex-tang-b}) in terms of
standard coordinates
$(x^i, y^i)$ are of the form
$$
\xi^C = y^k \frac{\partial \xi_{ij}}{\partial x^k}  dx^i\wedge dx^j +
\xi_{ij}dx^i \wedge dy^i,
\qquad v^C =  v^{ij} \frac{\partial}{\partial x^i}\wedge
\frac{\partial}{\partial y^j} +
y^k \frac{\partial v^{ij}}{\partial x^k}
\frac{\partial}{\partial y^i}\wedge
\frac{\partial}{\partial y^j}.
$$
  Consequently, the inner product of the lifts is
$$
i(v^C) \xi^C = v^{ij} \xi_{ij}.
$$
  But the lift of the inner product is
$$
(i(v)\xi)^C = y^k \frac{\partial}{\partial x^k} ( v^{ij} \xi_{ij}).
$$
  This may be explained by the fact that the complete lift of the tensor product
is not equal to the tensor product of complete lifts: $(u\otimes v)^C\ne
u^C\otimes v^C$ (cf. Remark~\ref{cl-prod} below).
\end{remark}

%
%

\subsection{The vertical lift of a tensor field}
\label{subsec-vl}

   The vertical lift of a tensor field
\begin{equation}
\label{u-tensor}
u=u_{i_1\dots i_k}^{j_1\dots j_\ell} dx^{i_1}\otimes \dots \otimes
dx^{i_k}\otimes \frac{\partial}{\partial x^{j_1}}
\otimes\dots\otimes \frac{\partial}{\partial x^{j_\ell}}
\end{equation}
in terms of the local coordinates
$(x^i, y^i)$ on~$TM$ is of the form~\cite{Yano-I}
$$
u^V=u_{i_1\dots i_k}^{j_1\dots j_\ell} dx^{i_1}\otimes \dots \otimes dx^{i_k}\otimes \frac{\partial}{\partial y^{j_1}}
\otimes\dots\otimes \frac{\partial}{\partial y^{j_\ell}}.
$$

  We generalize the notion of the vertical lift to the case of an arbitrary Frobenius
Weil algebra $\A$ in the following way.
  Consider a Jordan-H\"older basis  $\{e_0,\dots, e_n\}$ in $\A$.
  Let
$u\in {T}^{k,\ell}(M)$ be a  tensor field $(\ref{u-tensor})$ on $M$ and
$u^\A$ its analytic prolongation to $T^\A M$.
  We define the {\it vertical lift} $u^V\in {T}^{k,\ell}(T^\A M)$
of $u$ by
\begin{equation}
\label{vert-lift-def}
u^V=u^V_{\A}:=R(e_n u^\A).
\end{equation}
   The above definition does not depend on the choice of a Jordan-H\"older
basis because $e_n$ is a basis of the one-dimensional ideal $\Ao{}^h$ (see
Proposition~\ref{prop-fa}) and is fixed by the condition $p(e_n)=1$.

\begin{proposition}
  In terms of local coordinates $(x^{ia})$ on $T^\A M$,
\begin{equation}
\label{vert-lift}
u^V=u_{i_1\dots i_k}^{j_1\dots j_\ell} dx^{i_10}\otimes \dots \otimes
dx^{i_k0}\otimes \frac{\partial}{\partial x^{j_1n}}
\otimes\dots\otimes \frac{\partial}{\partial x^{j_\ell n}}.
\end{equation}
\end{proposition}

\begin{Proof}
   It follows from~(\ref{q-inv}) (see also  Remark~\ref{cl-inj}) that
$(e_n u^\A)_{i_1\dots i_k}^{j_1\dots j_\ell} e_{a_1}\dots e_{a_k} e^{b_1} \dots e^{b_\ell}$
vanishes for all values of indices
except for the case
$a_1=\ldots = a_k =0$ and
$b_1=\ldots = b_\ell =n$, when
it is equal to
$e_n (u_{i_1\dots i_k}^{j_1\dots j_\ell})^\A =
e_n u_{i_1\dots i_k}^{j_1\dots j_\ell}$.
  Contracting with~$p$, we obtain
$$
  (u^V)_{i_1a_1\dots i_ka_k}^{j_1b_1\dots j_\ell b_\ell}=\left\{
     \begin{array}{ll}
       u_{i_1\dots i_k}^{j_1\dots j_\ell},\quad &
       \hbox{if~}~a_1=\ldots = a_k =0,~~b_1=\ldots = b_\ell =n,\\[5pt]
       0,     & \hbox{otherwise}.
     \end{array}
  \right.
$$
\end{Proof}

  The following proposition follows immediately from~(\ref{vert-lift}).
\begin{proposition}
  For every $u,v\in {T}^{*,*}(M)$ we have
\begin{equation}
\label{vl-wedge}
 (u\otimes v)^V = u^V \otimes v^V.
\end{equation}
\end{proposition}

   It should be noted  that, in accordance with~(\ref{xiV}),
the vertical lift of a smooth function $f$ is equal to $f\circ
\pi_\A$~\cite{Shir12}.

\begin{remark}
\label{rem-vl}
For the discussion in Section~\ref{sect-wbpm} it is convenient to present here
a direct coordinate proof  of the fact that~(\ref{vert-lift}) defines a tensor
field on
$T^\A M$ (cf.~\cite{Vshjr-ljm}).
  Let $x^{i'}=x^{i'}(x^i)$ be a coordinate change on~$M$ and
$$
u_{i'_1\dots i'_k}^{j'_1\dots j'_\ell}
 = \frac{\partial x^{i_1}}{\partial x^{i'_1}} \dots
\frac{\partial x^{i_k}}{\partial x^{i'_k}}
\frac{\partial x^{j'_1}}{\partial x^{j_1}} \dots
\frac{\partial x^{j'_k}}{\partial x^{j_\ell}}
u_{i_1\dots i_k}^{j_1\dots j_\ell}.
$$
the corresponding transformation of components of a tensor $u$.
  It suffices to prove that
\begin{equation}
\label{dx-n}
\frac{\partial}{\partial x^{i'_sn}} =
\frac{\partial x^{i_s}}{\partial x^{i'_s}}
\frac{\partial}{\partial x^{i_sn}}.
\end{equation}
  Let us find the  change of coordinates
$x^{i'a}=x^{i'a}(x^{ib})$ on $T^\A M$ corresponding to a change $x^{i'}=x^{i'}(x^i)$.
  By (\ref{A-prol}), we have
\begin{equation}
\label{A-pr-c}
  X^{i'} = x^{i'0} + x^{i'\hat a} e_{\hat a} =
  x^{i'0} + \sum^h_{|p|=1}\frac{1}{p!}
    \frac{D^p x^{i'}}{Dx^p}\Xo{}^p.
\end{equation}
  Then, for $T^\A M$ we have
$$
  x^{i'0} = x^{i'0} (x^{i0}) ,\qquad
  x^{i'\hat a} = x^{i' \hat a} (x^{ib}), \quad \hat a=1,\dots, n.
$$
  Let us show that
$$
\frac{\partial x^{i'\hat a}}{\partial x^{ib}} =0 \quad \hbox{for} \quad
b> \hat a.
$$
  In fact, the coefficients  $\frac{1}{p!} \frac{D^p x^{i'}}{Dx^p}$
in (\ref{A-pr-c}) depend only on $x^{i}=x^{i0}$, while $x^{ib}$ occurs in an
expression
$(\Xo{}^1)^{p_1}\dots (\Xo{}^m)^{p_m}$  only as a coefficient of $e_b$
in $\Xo{}^i$.
  Since $\gamma_{c\hat d}^s=0$ for $c>s$,
the coefficient of $e_{\hat a}$ in the  expansion of
$(\Xo{}^1)^{p_1}\dots (\Xo{}^m)^{p_m}$ does not depend on~$x^{ib}$.
  Moreover, it can depend on $x^{i\hat a}$ only when $|p|=1$.
  The part of the sum (\ref{A-pr-c}) corresponding to $p=1$ is of
the form $\frac{\partial x^{i'}}{\partial x^{j}} x^{j\hat c} e_{\hat c}$, and
the variable $x^{i\hat a}$ appears in this expression only when
$\hat c=\hat a$.
  Therefore,
$$
\frac{\partial x^{i'\hat a}}{\partial x^{i\hat b}} =
\left\{
\begin{array}{ll}
0,\qquad &  \hbox{if} \quad
\hat b> \hat a,\\[5pt]
\displaystyle \frac{\partial x^{i'}}{\partial x^{i}}, &  \hbox{if} \quad
\hat b= \hat a,
\end{array}
\right.
$$
and the Jacobi matrix
$\left\| \displaystyle\frac{\partial x^{i'a}}{\partial x^{ib}} \right\|$
of the coordinate change
$x^{i'a}=x^{i'a}(x^{ib})$ on $T^\A M$ has the following block
structure:
\begin{equation}
\label{ch-coord}
\hbox{%
\begin{tabular}{|c|c|c|c|c|c|}
\hline
$\vphantom{\Biggl(} \frac{\partial x^{i'}}{\partial x^{i}} $ &
~$*$~ & ~$*$~ & \,\,$\cdots$\,\, & ~$*$~ & ~$*$~\\
\hline
~0~ & $\vphantom{\Biggl(} \frac{\partial x^{i'}}{\partial x^{i}} $ &
$*$ & $\cdots$ & $*$ & $*$\\
%
%
\hline
$0$ & $0$ & \,$\vphantom{\Biggl(}\ddots$\, &  $\ddots$ & $\vdots$ & $\vdots$\\
\hline
$\vdots$  & $\vphantom{\Biggl(}\ddots$  & $\ddots$  &$\ddots$  &
$\ddots$ & $\vdots$\\
\hline
$0$ &  $0$ & $\cdots$ & $0$ &
$\vphantom{\Biggl(}\frac{\partial x^{i'}}{\partial x^{i}} $ & $*$\\
\hline
$0$ & $0$ & $\cdots $ & $0$ & $0$ &
$\vphantom{\Biggl(}\frac{\partial x^{i'}}{\partial x^{i}} $ \\
\hline
\end{tabular}~~,}
\end{equation}
where $*$ denotes the blocks which are unessential for our
consideration.
  Now (\ref{dx-n}) is obvious.
\end{remark}

\begin{proposition}
\label{prop-vlbr}
  For  $u,v\in {\cal V}^*(M)$  we have
\begin{equation}
\label{vlbr}
\begin{array}{l}
a) \qquad [u,v]^V=[u^V, v^C] = [u^C, v^V];\\[5pt]
b) \qquad [u^V, v^V]=0.
\end{array}
\end{equation}
\end{proposition}

\begin{Proof}
   Let $u^\A$ and $v^\A$ be the $\A$-prolongations of
$u$ and $v$, respectively.
   Then
$$
[u,v]^V = R(e_n\, [u^\A,v^\A])= R([e_n u^\A, v^\A]) =
[ R(e_n u^\A), R(v^\A)] = [u^V, v^C].
$$
  Similarly, $[u,v]^V=[u^C, v^V]$.

  The second equality is proved in the same way:
$$
[u^V, v^V] = [R(e_n u^\A), R(e_n v^\A)] =
 R(e_n e_n \, [u^\A, v^\A]) = 0.
$$
\end{Proof}

  For the case of the tangent bundle~$TM$
the relations~(\ref{vlbr}) were proved by J. Grabowski and P.
Urba\'nski~\cite{G-U2, G-U3}.

\begin{remark}
\label{cl-prod}
  J. Grabowski and P. Urba\'nski have also proved~\cite{G-U2} that
for the tangent bundle one has
\begin{equation}
\label{uvC}
  (u\otimes v)^C = u^C \otimes v^V + u^V \otimes v^C
\end{equation}
for every $u,v\in {T}^{*,*}(M)$ (see also the monograph of
K.\,Yano and S.\,Ishihara~\cite{Yano-I}).
  This formula can be generalized for the case of
Frobenius Weil algebras of height $h>1$ in the following way.
  Let $u\in {T}^{k,\ell}(M)$ be a smooth tensor field and
$u^\A\in {T}^{k,\ell}_{\Adiff}(T^\A M)$  its analytic prolongation.
  Define the  {\it $a$-lift} $u^{(a)}\in {T}^{k,\ell}(T^\A M)$
of $u$ for a fixed basis $\{e_a\}$ by
\begin{equation}
\label{a-lift}
u^{(a)}:=R(e_a u^\A).
\end{equation}
   In particular, when  $\{e_a\}$ is the Jordan-H\"older basis,
$u^{(0)}=u^C$, $u^{(n)}=u^V$.

   For  $u,v\in {T}^{*,*}(M)$
the following generalization of the relation~(\ref{uvC}) holds:
\begin{equation}
\label{cl-prod-f}
  (u\otimes v)^C = \sum_{a,b} q^{ab} u^{(a)} \otimes v^{(b)}.
\end{equation}

  In fact, let $u^\A, v^\A\in {T}^{*,*}_{\Adiff}(T^\A M)$
be the analytic prolongations of~$u$ and~$v$, respectively.
  To prove~(\ref{cl-prod-f}) it suffices to verify that
$$
p(XY) = q^{ab} p(e_a X) p(e_b Y).
$$
  We have
\begin{gather*}
p(XY) = X^c Y^d p(e_c e_d) = X^c Y^d q_{cd},
\\
q^{ab} p(e_a X) p(e_b Y) = q^{ab} X^c p(e_a e_c) Y^d p(e_b e_d) = q^{ab} X^c
Y^d q_{ac} q_{bd} = X^c Y^d \delta^b_c q_{bd} = X^c Y^d q_{cd}.
\end{gather*}

  In the similar manner one can prove that
\begin{equation}
\label{cl-prod-f2}
  (u\otimes v)^{(a)} = \sum_{b,d}\gamma^{bd}_a
          u^{(b)} \otimes v^{(d)}.
  \end{equation}

  The tangent bundle $T^n(M)= T^{\R(\varepsilon^n)}M$ of  order $n$ is equivalent
to the Weil bundle corresponding to  the algebra of plural numbers
$\R(\varepsilon^n)$.
  In this case the so called $\lambda$-lifts
$u^{(\lambda)} = R(\varepsilon^\lambda u^{\R(\varepsilon^n)})$
of a tensor field $u$ on $M$ ($\lambda=0,\dots, n$) are defined.
   These lifts were considered, for example,
in the papers of Ch.-S.~Houh and S.~Ishihara~\cite{H-I} and
A.P.~Shirokov~\cite{Shir5}.
   Relation~(\ref{cl-prod-f2}), in this case takes the form~\cite{H-I, Shir5}
$$
  (u\otimes v)^{(\lambda)} = \sum_{\kappa=0}^\lambda
   u^{(\kappa)}\otimes v^{(\lambda-\kappa)}.
$$
\end{remark}

\begin{proposition}
\label{prop-lie-der}
  Let $v\in {\cal V}^1(M)$  be a vector field on $M$.
  The Lie derivative with respect to $v$ has the following properties:
\begin{equation}
\label{lie-lift}
\begin{array}{l}
a) \qquad ({\cal L}_v t)^C = {\cal L}_{v^C} t^C
; \\[5pt]
b) \qquad ({\cal L}_v t)^V = {\cal L}_{v^C} t^V= {\cal L}_{v^V} t^C
;\\[5pt]
c) \qquad {\cal L}_{v^V} t^V = 0
\end{array}
\end{equation}
for every $t\in {T}^{*,*}(M)$.
\end{proposition}

\begin{Proof}
a)  Recall that the components of the Lie derivative of a tensor field
$$
t= t_{i_1\dots i_k}^{j_1\dots j_\ell} dx^{i_1}\otimes \dots
\otimes dx^{i_k}\otimes \frac{\partial}{\partial x^{j_1}}
\otimes\dots\otimes \frac{\partial}{\partial x^{j_\ell}}
$$
with respect to a vector field $v$ are of the form~\cite{Post, Yano}
\begin{multline}
\label{lie-der}
\displaystyle
({\cal L}_v t)_{i_1\dots i_k}^{j_1\dots j_\ell} =  \frac{\partial t_{i_1\dots
i_k}^{j_1\dots j_\ell}}{\partial x^m} v^m + t_{mi_2\dots i_k}^{j_1\dots j_\ell}
\frac{\partial v^m}{\partial x^{i_1}} + \ldots +
t_{i_1\dots i_{k-1}m}^{j_1\dots j_\ell} \frac{\partial v^m}{\partial x^{i_k}}\,  \\[12pt]
\displaystyle
-\, t_{i_1\dots i_k}^{mj_2\dots j_\ell} \frac{\partial v^{j_1}}{\partial x^{m}}
- \ldots -
 t_{i_1\dots i_k}^{j_1\dots j_{\ell-1}m} \frac{\partial v^{j_\ell}}{\partial x^{m}}.
\end{multline}

  For simplicity, we  prove~(\ref{lie-lift}) in  the case $k=\ell=1$.
  In the general case this formula is proved in the same way.
  Let $t=t_i^j dx^{i}\otimes \frac{\partial}{\partial x^{j}}$ be a
 $(1,1)$-tensor field on $M$ and let $T$ and $V$ be the analytic
prolongations of $t$ and $v$, respectively.
  Then the analytic prolongations satisfy the rela\-tion
\begin{equation}
\label{lie11}
\displaystyle
({\cal L}_V T)_{i}^{j}
= \frac{\partial T_{i}^{j}}{\partial X^m} V^m +
T_{m}^{j} \frac{\partial V^m}{\partial X^{i}}
- T_{i}^{m} \frac{\partial V^{j}}{\partial X^{m}}.
\end{equation}
  Recall that if $V^i= v^{ia} e_a = v^i_b e^b$, then
$(v^C)^{ia}= v^{ia}=v^i_b q^{ab}$ by (\ref{real-vector}).

  Let us multiply each side of~(\ref{lie11}) by $e_a e^b$ and
then contract the result with $p$.
  On the left-hand side we have $(({\cal L}_v t)^C)_{ia}^{jb}$, and
it remains to prove that  the right-hand is of the form
$$
\frac{\partial t_{ia}^{jb}}{\partial x^{mc}} v^{mc} +
t_{mc}^{jb} \frac{\partial v^{mc}}{\partial x^{ia}}
- t_{ia}^{mc} \frac{\partial v^{jb}}{\partial x^{mc}}
$$
where
$t_{ia}^{jb}=(t^C)_{ia}^{jb}$ and $v^{ia}=(v^C)^{ia}$.

  First, let $T_i^j e_a e^b = (T_{ia}^{jb})_c e^c$.
  Then, by virtue of~(\ref{der-dir}), (\ref{partialF} and~(\ref{pee-up}),
$$
\displaystyle
p \Bigl( \frac{\partial T_{i}^{j}}{\partial X^m} V^m e_a e^b \Bigr) =
p \Bigl( \frac{\partial (T_{ia}^{jb})_c}{\partial x^{md}}
\delta^d e^c v^m_s e^s \Bigr)
=
\frac{\partial (T_{ia}^{jb})_d}{\partial x^{mc}}
\delta^d  v^m_s q^{cs}
= \frac{\partial t_{ia}^{jb}}{\partial x^{mc}} v^{mc}.
$$

 Second, by~(\ref{u-s3}) we have
$$
\begin{array}{l}
\displaystyle
T_{m}^{j} \frac{\partial V^m}{\partial X^{i}} e_a e^b = T_{m}^{j} \delta^d
\frac{\partial v^{mg}}{\partial x^{id}} e_g e_a e^b = T_{m}^{j} \delta^d
\frac{\partial v^{mg}}{\partial x^{id}} \gamma_{ag}^c e_c e^b =
 T_{m}^{j} e_c e^b  \frac{\partial v^{mc}}{\partial x^{ia}} = (T_{mc}^{jb})_d
e^d \frac{\partial v^{mc}}{\partial x^{ia}} .
\end{array}
$$
  Whence
$$
p \Bigl(  T_{m}^{j} \frac{\partial V^m}{\partial X^{i}} e_a e^b \Bigr)=
t_{mc}^{jb} \frac{\partial v^{mc}}{\partial x^{ia}} .
$$

  Finally, let $V^j e^b = V^{jb}_c e^c$.
  Then
$$
\begin{array}{l}
\displaystyle
T_i^m \frac{\partial V^j}{\partial X^m} e_a e^b =
T_i^m e_a \delta^d \frac{\partial V^{jb}_c}{\partial x^{md}}  e^c =
T_i^m e_a e^c \delta^d \frac{\partial V^{jb}_d}{\partial x^{mc}} =
T_i^m e_a e^c \frac{\partial v^{jb}}{\partial x^{mc}}.
\end{array}
$$
  Therefore,
$$
p \Bigl( T_i^m \frac{\partial V^j}{\partial X^m} e_a e^b \Bigr)=
t_{ia}^{mc} \frac{\partial v^{jb}}{\partial x^{mc}} .
$$

  Formulas (\ref{lie-lift}) b) and~c) are  proved in the same way.
\end{Proof}

  For   tangent bundles  relations (\ref{lie-lift}) has
been  proved in~\cite{Yano-I} and~\cite{Vai4}.

%
%
%

\section{Weil bundles of Poisson manifolds}
\label{sect-wbpm}

   In this section we consider the
complete and the vertical lifts of a Poisson tensor $w$ and establish some
properties of the  Poisson structures arised. In particular, we compute the
modular classes of these structures.

\subsection{The complete lift of a Poisson tensor}
\label{subsec-cl-ps}

  Let $(M,w)$ be a Poisson manifold, and let
$w^C$ be the complete lift of $w$ to $T^\A M$.
  By virtue of~(\ref{ww0}) and~(\ref{lift-Sch-br}),
$$
  [w^C,w^C]=0.
$$
Hence $w^C$ is a Poisson tensor on $T^\A M$.

   Let $w^\A$ be the analytic prolongation of~$w$
and let $(w^{ij})^\A = (w^{ij})^s e_s$ be the expansions in terms of a basis in
$\A$.
   Then $(w^{ij})^\A e^a e^b = (w^{ij})^s e_s e^a e^b $, therefore
$(w^C)^{iajb}=(w^{ij})^s  \gamma^{ab}_s$ by~(\ref{p-gamma}).
   Thus, the components of $w^C$ in terms of the local
coordinates on $T^\A $M are as follows
\begin{equation}
\label{wC}
  (w^C)^{iajb} = (w^{ij})^s  \gamma^{ab}_s .
%
%
\end{equation}

\begin{proposition}
\label{TA-sym}
  If $(M,w)$ is a symplectic manifold, then
$(T^\A M, w^C)$ also is a symplectic manifold.
\end{proposition}

\begin{Proof}
  By the Darboux theorem~\cite{daS-W}, we can choose a local coordinate system $(x^i)$
on~$M$, in terms of which the components $w^{ij}$ are constant.
  Then, by (\ref{A-prol}), the analytic prolongations
$(w^{ij})^\A$ coincide with~$w^{ij}$ and so~$(w^{ij})^0=w^{ij}$ and~$(w^{ij})^k=0$ for~$k\ge 1$.
  Therefore $w^{iajb} = w^{ij} \gamma^{ab}_0 = w^{ij}
q^{ab}$.
  Thus,
$$
  w^C = w^{ij} q^{ab} \frac{\partial}{\partial x^{ia}}\wedge
   \frac{\partial}{\partial x^{jb}}
$$
where the matrix $\|(w^C)^{iajb}\|=\|w^{ij} q^{ab}\|$ is nondegenerate as
the tensor (Kronecker) product of nondegenerate matrices
$\|w^{ij}\|$ and $\| q^{ab}\|$.
\end{Proof}

  One can easily see that if~$\omega=\omega_{ij}\,dx^i\wedge dx^j$ is the
symplectic form on $M$ corresponding to $w$, then the complete lift
$$
  \omega^C = \omega_{ij} q_{ab} \, dx^{ia}\wedge dx^{jb}
$$
is the symplectic form on~$T^\A M$ corresponding to~$w^C$.
 In the same way one can prove that if $(M,w)$ is a regular Poisson manifold, then
$(T^\A M, w^C)$ is  also a regular Poisson manifold.

\begin{remark}
   The fact that for a Frobenius Weil algebra $(\A,q)$ the total space $T^\A M$ of
 a symplectic manifold
$(M,w)$  carries a natural symplectic structure,
was pointed out by A.V.~Brailov~\cite{Bra}.
\end{remark}

%
%
\begin{example}
   Let $(M,w)$ be a Poisson manifold and  $TM$  its tangent bundle.
      The complete lift of $w$ corresponding to
$p_{(0)}$ (Example~\ref{ex-tang-b}) is of the form
\begin{equation}
\label{wc-tan}
\displaystyle
  w^C= w^{ij} \frac{\partial }{\partial x^i}
   \wedge \frac{\partial }{\partial y^j} +
     y^k \frac{\partial w^{ij}}{\partial x^k}
     \frac{\partial }{\partial y^i} \wedge
      \frac{\partial }{\partial y^j},
\end{equation}
and the complete lift $w^C$, corresponding to  $p_{(1)}$, is
$$
\displaystyle
  w^C= w^{ij} \frac{\partial }{\partial x^i}
   \wedge \frac{\partial }{\partial y^j} +
    \biggl( y^k \frac{\partial w^{ij}}{\partial x^k} -
       w^{ij}\biggr)
     \frac{\partial}{\partial y^i} \wedge
      \frac{\partial}{\partial y^j}.
$$
  Poisson bivector~(\ref{wc-tan}) was studied by several authors, e.g.,
T.J.~Courant~\cite{Cou}, J.~Gra\-bowski and P.~Urba\'nski~\cite{G-U, G-U2,
G-U3}, G.~Mitric and I.~Vaisman~\cite{M-V, Vai4}.
\end{example}

\begin{proposition}
Let $w$ be a bivector field on $M$.
   The complete lift $w^C$  is a Poisson bivector on $T^\A M$
if and only if $w$ is a Poisson bivector on $M$.
\end{proposition}

\begin{Proof}
   Follows from~(\ref{lift-Sch-br}) and Remark~\ref{cl-inj}.
\end{Proof}

\begin{proposition}
\label{prop-cl-pmap}
Let $(M,w)$ and $(M',w')$ be  Poisson manifolds and
  let $\varphi : (M,w) \to (M',w')$ be a Poisson map.
  Then $T^\A\varphi : (T^\A M,w^C) \to (T^\A M',(w')^C)$ is also  a Poisson map.
\end{proposition}

\begin{Proof}
  Follows from (\ref{pmap}) and Proposition~\ref{prop-real-contrav}.
\end{Proof}

\begin{proposition}
  The complete lift of multivector fields induces the homomorphism of
Poisson cohomology spaces
\begin{equation}
\label{wc-hom}
  H_{P}^*(M,w)\, \longrightarrow \, H_{P}^* (T^\A M,w^C),
  \qquad   [u] \longmapsto [u^C].
\end{equation}
\end{proposition}

\begin{Proof}
  From (\ref{lift-Sch-br}) it follows that
\begin{equation}
\label{sigmaw-c}
  (\sigma_w u)^C = \sigma_{w^C}\, u^C,
\end{equation}
which implies that the map (\ref{wc-hom}) is well-defined.
\end{Proof}

  It follows from  Proposition~\ref{prop-Sch-br} that if $f\in C^\infty(M)$
is a Casimir function of $w$ then $f^C=R(f^\A)$ and $f^V=f\circ \pi_\A$ are
Casimir functions of $w^C$.
  For an arbitrary smooth function $f$ on $M$ one has
$$
  (X^w_f)^C=X^{w^C}_{f^C}.
$$

\begin{proposition}
\label{cl0}
  If $p(1_\A)\ne0$, then the homomorphism {\rm (\ref{wc-hom})}  in the dimension $0$
\begin{equation}
\label{pc0}
H_{P}^0(M,w)\, \longrightarrow \,   H_{P}^0 (T^\A M,w^C)
\end{equation}
is a monomorphism.
  If $p(1_\A)=0$ then the kernel of~{\rm(\ref{pc0})} is the space of constant functions.
\end{proposition}

\begin{Proof}
  Follows from  Remark~\ref{cl-inj-cov}.
\end{Proof}

\begin{proposition}
\label{pm-cohom-cd}
  Let $(M,w)$ be a Poisson manifold.
  For every Frobenius Weil al\-gebra~$(\A,q)$ the following diagram is commutative
$$
\xymatrix@R=15mm{%
    {{\cal V}^*(M)} \ar[rr]^-{C} & &
         {{\cal V}^*(T^\A M)}  \\
    {\Omega^*(M)} \ar[rr]^-{C}\ar[u]^{\widetilde w} &
     & {\,\,\Omega^*(T^\A M)~} \ar[u]_{\widetilde{w^C}}, \\
}
$$
where the horizontal arrows mean the complete lift.
\end{proposition}

\begin{Proof}
    Define $\gamma^{a_1\dots a_k}_b$ by
$$
e^{a_1}\dots e^{a_k} = \gamma^{a_1\dots a_k}_b e^b.
$$

  Consider an exterior form $\xi\in \Omega^k(M)$ and a multivector field
$v\in{\cal V}^k(M)$.
  In terms of local coordinates,
$$
\xi=\xi_{i_1\dots i_k} dx^{i_1}\wedge\ldots\wedge dx^{i_k}
\quad \hbox{and} \quad
v=v^{j_1\dots j_k} \frac{\partial}{\partial x^{j_1}}
\wedge\ldots\wedge \frac{\partial}{\partial x^{j_k}}.
$$
  Let $(\xi_{i_1\dots i_k})^\A = \xi^f_{i_1\dots i_k} e_f$
and $(v^{j_1\dots j_k})^\A = (v^{j_1\dots j_k})^s e_s$ be the analytic prolongations
of $\xi_{i_1\dots i_k}$ and
$v^{j_1\dots j_k}$, respectively.
  Then
$$
(\xi^C)_{i_1a_1\dots i_ka_k}
= p(\xi^f_{i_1\dots i_k} e_f e_{a_1}\dots e_{a_k})
= \xi^f_{i_1\dots i_k} \gamma_{f a_1\dots a_k}^d p_d,
$$
and
\begin{multline*}
(v^C)^{j_1 b_1\dots j_k b_k} \\= p((v^{j_1\dots j_k})^s e_s e^{b_1}\dots
e^{b_k})
 = (v^{j_1\dots j_k})^s p( q_{sf} e^f e^{b_1}\dots e^{b_k}) =
(v^{j_1\dots j_k})^s q_{sf} \gamma^{f {b_1}\dots {b_k}}_d \delta^d,
\end{multline*}
where $\gamma_{a_1\dots a_k}^b$ are defined by (\ref{gammaup}).

  In our case $v=\widetilde w\xi$.
  According to (\ref{tw-f}),
$v^{j_1\dots j_k}= (-1)^k w^{i_1j_1}\dots w^{i_kj_k} \theta_{i_1\dots i_k}$.
  Thus, we need to show that
$$
  (w^C)^{i_1a_1j_1b_1}  \dots   (w^C)^{i_ka_kj_kb_k}
  (\xi^C)_{i_1a_1\dots i_ka_k} =
  (v^C)^{j_1 b_1\dots j_k b_k}.
$$
  We have
$$
  (w^C)^{i_1a_1j_1b_1}  \dots   (w^C)^{i_ka_kj_kb_k}
  (\xi^C)_{i_1a_1\dots i_ka_k} =
  (w^{i_1j_1})^{c_1} \gamma^{a_1b_1}_{c_1} \dots   (w^{i_kj_k})^{c_k}
  \gamma^{a_kb_k}_{c_k}
   \xi^f_{i_1\dots i_k} \gamma_{f a_1\dots a_k}^d p_d.
$$
  On the other hand, it follows from  Proposition~\ref{A-prol-prop} that
$$
(v^{j_1\dots j_k})^\A =
  (w^{i_1j_1})^{c_1} \dots (w^{i_kj_k})^{c_k}
  \xi^f_{i_1\dots i_k} e_f e_{c_1}\dots e_{c_k}.
$$
Therefore
$$
(v^C)^{j_1b_1\dots j_kb_k} =
  (w^{i_1j_1})^{c_1} \dots (w^{i_kj_k})^{c_k}
  \xi^f_{i_1\dots i_k} p(e_f e_{c_1}\dots e_{c_k}
   e^{b_1}\dots e^{b_k}).
$$
   Thus, it suffices to prove that
\begin{equation}
\label{gggp}
  \gamma^{a_1b_1}_{c_1} \dots \gamma^{a_kb_k}_{c_k}
   \gamma_{f a_1\dots a_k}^d p_d =
   p(e_f e_{c_1}\dots e_{c_k}
   e^{b_1}\dots e^{b_k}).
\end{equation}
  We have
\begin{multline*}
p(e_f e_{c_1}\dots e_{c_k} e^{b_1}\dots e^{b_k}) \\= q^{a_1b_1} \dots
q^{a_kb_k} p(e_{a_1}\dots e_{a_k} e_f e_{c_1}\dots e_{c_k}) = q^{a_1b_1} \dots
q^{a_kb_k}
\gamma_{f{a_1}\dots a_k}^d \gamma_{d{c_1}\dots c_k}^s p_s.
\end{multline*}
  According to~(\ref{gamma2-delta}),
$q^{a_1b_1}=\gamma^{a_1b_1}_{c_1} \delta^{c_1}$, \dots,
$q^{a_kb_k}=\gamma^{a_kb_k}_{c_k} \delta^{c_k}$.
   Since for a Jordan-H\"older basis $\delta^0=1$ and $\delta^{\hat a}=0$, it follows that
$\gamma_{d{c_1}\dots c_k}^s p_s \delta^{c_1}\dots
\delta^{c_k} = \delta^s_d p_s = p_d$.
  Hence the result.
\end{Proof}

\begin{theorem}
    The following diagram of morphisms of complexes is commutative
\begin{equation}
\label{cl-mor}
\xymatrix@R=15mm{%
    {({\cal V}^*(M), \sigma_w)} \ar[rr]^-{C} &
    &     {({\cal V}^*(T^\A M), \sigma_{w^C})}\\
    {(\Omega^*(M),d)} \ar[rr]^-{C}\ar[u]^{\widetilde w} &
    &  {(\Omega^*(T^\A M),d)} \ar[u]_{\widetilde{w^C}} \\
}
\end{equation}
\end{theorem}

\begin{Proof}
   The commutativity of the three-dimensional diagram
$$
\xymatrix{%
 & {{\cal V}^{k}(M)} \ar[rr]^{C}
& & {{\cal V}^k(T^\A M)} \\
{{\cal V}^{k-1}(M)} \ar[ur]^{\sigma_w} \ar[rr]^(0.6){C}
 & \ar[r] \ar[u] & {{\cal V}^{k-1}(T^\A M)} \ar[ur]_{\sigma_{w^C}}\\
 & \Omega^k(M) \ar@{-}[r]^(0.7){C}
\ar@{-}[u]^(0.8){\widetilde w} & \ar[r] & \Omega^k(T^\A M)
\ar[uu]_{\widetilde {w^C}} \\
\Omega^{k-1}(M)\ar[uu]^{\widetilde w} \ar[rr]^{C}
\ar[ur]^d & & \Omega^{k-1}(T^\A M) \ar[ur]_d
\ar[uu]_(0.35){\widetilde {w^C}}\\
}
$$
follows from~(\ref{lift-ext-dif}), (\ref{sigmaw-c})
and~(\ref{d-sigma}).
\end{Proof}

\begin{corollary}
\label{cl-symp-p}
  Let $(\A,q)$ be a Frobenius Weil algebra and  let $p$ be the corresponding
Frobenius covector.
  Let $(M, w)$ be a symplectic manifold.
  If $p(1_\A)\ne 0$, then~{\rm(\ref{wc-hom})} is an isomorphism.
  If $p(1_\A)=0$, then~{\rm(\ref{wc-hom})} is the zero map.
\end{corollary}

\begin{Proof}
   For a symplectic manifold the maps $\widetilde w$ and
$\widetilde{w^C}$ are isomorphisms.
  Thus, the vertical arrows in~(\ref{cl-mor})  are isomorphisms.
  The rest of the proof follows from Theorem~\ref{cl-dR}.
\end{Proof}

\begin{remark}
   It has been  shown in \cite{M-V} and \cite{Vai3} that for the tangent bundle
$TM$ of a Poisson manifold $(M,w)$  the complete lift
$w^C$~(\ref{wc-tan}) is an exact Poisson structure.
   In fact, for $E=y^i \frac{\partial}{\partial y^i}$ we have
$w^C=\sigma_{w^C} E = [w^C, E]$.
   One can easily verify that the vertical lift
$w^V= w^{ij} \frac{\partial}{\partial y^i}
\wedge \frac{\partial }{\partial y^j}$ of $w$ is also exact:
$w^V = [w^V, E]$.
\end{remark}

  In the case of an arbitrary Frobenius Weil algebra,
 Corollary~\ref{cl-symp-p} implies the following proposition.

\begin{proposition}
   Let $(\A,q)$ be a Frobenius Weil algebra, and let
$p$ be the corresponding Frobenius covector such that $p(1_\A)= 0$.
   If $(M,w)$ is a symplectic manifold then
$(T^\A M, w^C)$ is an exact symplectic manifold.
\end{proposition}


   The following example shows that homomorphism~{\rm(\ref{wc-hom})}
may be a monomorphism or have a nonzero kernel depending on the dimension of
the cohomology space.

\begin{example}
\label{ex-torus-Rk}
   Let $\T^2=\Sb^1\times \Sb^1$ be the two-dimensional torus and let $M=\T^2\times \R^k$.
   We denote the standard angle coordinates on the torus by $(x^1, x^2)$ and
the standard coordinates on $\R^k$ by $(t^1,\dots, t^k)$.
   Consider the case of the algebra $\R(\varepsilon)$ with  Frobenius covector
$p_{(0)}$ (see Example~\ref{ex-tang-b}).
  The  corresponding Weil bundle of $M$ coincides with its tangent bundle: $T^{\R(\varepsilon)}M=TM$.
   The tangent bundle of $M=\T^2\times \R^k$ is trivial:
$TM\cong M \times (\R^2\times \R^{k})$.
   Denote by $(y^1, y^2)$ and $(s^1,\dots, s^k)$ the standard coordinates in $\R^2$ and
$\R^k$,  respectively.

   The bivector field
$w=\frac{\partial}{\partial x^1}\wedge \frac{\partial}{\partial x^2}$
defines a regular Poisson structure on~$M$.
   The complete lift of $w$ is of the form $w^C=\frac{\partial}{\partial x^1}\wedge
\frac{\partial}{\partial y^2}
+ \frac{\partial}{\partial y^1}\wedge \frac{\partial}{\partial x^2}$.

   In the dimension 0, by~(\ref{reg-pr}) we have
$$
\begin{array}{l}
H^0_P(M,w)\cong C^\infty(\R^k)\cong \R\oplus
C^\infty_0(\R^k),\\[5pt]
H^0_P(TM,w^C)\cong C^\infty(\R^{2k}),
\end{array}
$$
where  $C^\infty_0(\R^k)$ is the subring of $C^\infty(\R^k)$, consisting of
smooth functions vanishing at
$0\in \R^k$.
   The complete lift of a function $f\in C^\infty(M)$ is
$f^C = y^i \frac{\partial f}{\partial x^i}+
s^a \frac{\partial f}{\partial t^a}$.
   The kernel of homomorphism~(\ref{wc-hom})
consists of constant functions (see Remark~\ref{cl-inj-cov}).
   Hence, its image is isomorphic to $C^\infty_0(\R^k)$.

   In the dimension 1, we have
$$
\begin{array}{l}
H^1_P(M,w)\cong
{\cal V}^1(\R^k) \oplus C^\infty(\R^k)\oplus  C^\infty(\R^k),\\[5pt]
H^1_P(TM,w^C)\cong
{\cal V}^1(\R^{2k}) \oplus C^\infty(\R^{2k})\oplus  C^\infty(\R^{2k}).
\end{array}
$$
  The cohomology classes of the following vector fields
form a complete system of  generators of  $H^1_P(M,w)$:
$$
f^i \frac{\partial}{\partial x^i}, \quad f^a\frac{\partial}{\partial t^a},
$$
where $f^i, f^a\in C^\infty(\R^k)$.
   The complete lifts of the above indicated vector fields, respecti\-vely,  are
$$
  \Bigl(f^i \frac{\partial}{\partial x^i}\Bigr)^C=
f^i \frac{\partial}{\partial x^i}
+ s^b \frac{\partial f^i}{\partial t^b}
\frac{\partial}{\partial y^i},
\qquad
\Bigl(f^a\frac{\partial}{\partial t^a}\Bigr)^C =
f^a\frac{\partial}{\partial t^a}
+ s^b \frac{\partial f^a}{\partial t^b}
\frac{\partial}{\partial s^a}.
$$
   One can easily see that
$f^1 \frac{\partial}{\partial x^1}= -\sigma_{w^C} (y^2f^1)$,
$f^2 \frac{\partial}{\partial x^2}= \sigma_{w^C} (y^1f^2)$,
and that the classes
$[s^b \frac{\partial f^i}{\partial t^b}
\frac{\partial}{\partial y^i}]$ (where at least one of $f^i$ is not a constant function) and
$[(f^a\frac{\partial}{\partial t^a})^C]$ are nonzero classes in $H^1_P(TM, w^C)$.
   Therefore, in the dimension 1, the kernel of homomorphism~(\ref{wc-hom})
is isomorphic to $\R\oplus \R$, and the image of~(\ref{wc-hom}) is isomorphic
to
${\cal V}^1(\R^k)\oplus C^\infty_0(\R^k)\oplus C^\infty_0(\R^k)$.

   In the dimension 2, we have
$$
\begin{array}{l}
H^2_P(M,w)\cong
{\cal V}^2(\R^k)\oplus {\cal V}^1(\R^k)\oplus
{\cal V}^1(\R^k)\oplus C^\infty(\R^k),\\[5pt]
H^2_P(TM,w^C)\cong
{\cal V}^2(\R^{2k})\oplus {\cal V}^1(\R^{2k})\oplus
{\cal V}^1(\R^{2k})\oplus C^\infty(\R^{2k}).
\end{array}
$$
  The cohomology classes of the following bivector
fields are generators of $H^2_P(M,w)$:
$$
f \frac{\partial}{\partial x^1} \wedge \frac{\partial}{\partial x^2},
\quad
f^{ia} \frac{\partial}{\partial x^i}\wedge
\frac{\partial}{\partial t^a}\quad \hbox{and}
\quad
f^{ab}\frac{\partial}{\partial t^a}\wedge
\frac{\partial}{\partial t^b},
$$
where $f, f^{ia}, f^{ab}\in C^\infty(\R^k)$.
   Among the cohomology classes
defined by these  bivector fields,  only the classes
 $[(f \frac{\partial}{\partial x^1}
\wedge \frac{\partial}{\partial x^2})^C]$ with
$f={\rm const}$ are zero.
   Thus, in the dimension 2, the kernel of~(\ref{wc-hom})
is isomorphic to $\R$, and the image of~(\ref{wc-hom}) is isomorphic to
$$
{\cal V}^2(\R^k)\oplus {\cal V}^1(\R^k)\oplus {\cal V}^1(\R^k)
\oplus C^\infty_0(\R^k).
$$

   In the dimensions $s\ge 3$, for every generator of $H^s_P(M,w)$
one can find a repre\-sen\-tative which has  the form
$u=u^a\wedge\frac{\partial}{\partial t^a}$, $u^a\in{\cal V}^{s-1}(M)$.
   Therefore, the cohomology class of $u^C$ is nonzero.
   Thus, in each dimension $s=3, \dots, k+2$, the complete lift
induces a monomorphism of Poisson cohomology spaces.
\end{example}

\subsection{The vertical lift of a Poisson tensor}
\label{subsec-vl-ps}

   It follows from~(\ref{vlbr}) that the vertical lift~$w^V$
of a Poisson tensor (as well as of any  bivector) $w$  is a Poisson tensor
on~$T^\A M$.

  The following example shows that the
cohomology of $\sigma_{w^V}$ depends on the choice of $w\in {\cal V}^2(M)$.

\begin{example}
  Let $M=\R^{2m}\times \R^{2k}$ and let $w_1$ be the regular
 structure induced by the standard symplectic structure on~$\R^{2m}$,
$w_2$  the regular Poisson structure induced by
the standard symplectic structure on~$\R^{2k}$. The sum
$w_1+w_2$ is the standard symplectic structure on~$M$.

   In the  case of the algebra $\R(\varepsilon)$,
$TM\cong \R^{4(m+k)}=\R^{2m+2k}\times \R^{2m}\times \R^{2k}$.
   By virtue of~(\ref{reg-pr}), one can obtain
$$
\begin{array}{l}
  H^r_P(TM, w_1^V) \cong {\cal V}^r (\R^{2m+2k}\times \R^{2k}),\\[5pt]
  H^r_P(TM, w_2^V) \cong {\cal V}^r (\R^{2m+2k}\times \R^{2m}),\\[5pt]
  H^r_P(TM, (w_1+w_2)^V) \cong {\cal V}^r (\R^{2m+2k}).
\end{array}
$$
   If $k>m$, then $H^r_P(TM, (w_1+w_2)^V)=0$ for $r>2m+2k$,
$H^r_P(TM, w_2^V)=0$ for $r>4m+2k$, and
$H^r_P(TM, w_1^V)=0$ for $r>4k+2m$.

  In  the case of an arbitrary Frobenius Weil algebra $\A$
of dimension~$n+1$,
$$
\begin{array}{l}
  H^r_P(T^\A M, w_1^V) \cong {\cal V}^r (\R^{2(m+k)n}\times \R^{2k}),\\[5pt]
  H^r_P(T^\A M, w_2^V) \cong {\cal V}^r (\R^{2(m+k)n}\times \R^{2m}),\\[5pt]
  H^r_P(T^\A M, (w_1+w_2)^V) \cong {\cal V}^r (\R^{2(m+k)n}).
\end{array}
$$
\end{example}

\begin{proposition}
\label{prop-vl-pmap}
Let $(M,w)$ and $(M',w')$ be  Poisson manifolds and
  let $\varphi : (M,w) \to (M',w')$  be a Poisson map.
  Then $T^\A\varphi : (T^\A M,w^V) \to (T^\A M',(w')^V)$ is also  a Poisson map.
\end{proposition}

\begin{Proof}
  Similar to that of Proposition~\ref{prop-cl-pmap}.
\end{Proof}

  Obviously, if $(M,w)$ is a regular Poisson manifold, then
$(T^\A M, w^V)$ is also a regular Poisson manifold.

\begin{proposition}
  The vertical lift of multivector fields induces the homomorphism of
Poisson cohomology spaces
\begin{equation}
\label{wv-hom}
  H_{P}^*(M,w)\, \longrightarrow \,
  H_{P}^* (T^\A M,w^V), \qquad
  [u] \longmapsto [u^V].
\end{equation}
\end{proposition}

\begin{Proof}
  From (\ref{vlbr}) it follows that if
$\sigma_w u=[w,u]=0$, then $\sigma_{w^V} u^V=[w^V, u^V]=0$, and if
$u =\sigma_w v= [w,v]$, then
$u^V  = \sigma_{w^V} v^C= [w^V, v^C]$.
\end{Proof}

  If $f\in C^\infty(M)$
is a Casimir function of $w$, then $f^C$ and $f^V$ are Casimir functions of
$w^V$.
  For an arbitrary smooth function $f$ on $M$, it follows from
Proposition~\ref{prop-vlbr} that
$$
  (X^w_f)^V=X^{w^V}_{f^C}  = X^{w^C}_{f^V}.
$$

\begin{proposition}
\label{vl0}
 In the dimension $0$,  homomorphism $(\ref{wv-hom})$
$$
H_{P}^0(M,w)\, \longrightarrow \,   H_{P}^0 (T^\A M,w^V)
$$
is a monomorphism.
\end{proposition}

\begin{Proof}
  If $[f^V]=[f\circ \pi_\A]=0$, then
   $f=0$.
\end{Proof}

   The next example shows that the homomorphism~(\ref{wv-hom}) may be the
zero map in every dimension except for zero dimension.

\begin{example}
   Let $(M,w)$ be a symplectic manifold and
let $(x^i)$ be the local coordinate system on $M$ with respect to which the
components $w^{ij}$ are constant.
  Let
$(x^i, y^i)$ be the induced local coordinate system on~$TM$.

   The space $H_P^0(M,w)$ coincides with the set of  constant functions on~$M$.
   The space $H^0_P(TM, w^V)$ coincides with the set of  constant
functions on~$TM$ having constant values along fibers.
  Therefore, it is isomorphic
to~$C^\infty(M)$.
  By Proposition~\ref{vl0}, the vertical lift induces a monomorphism
$H_P^0(M,w) \cong \R \to H^0_P(TM, w^V) \cong C^\infty(M)$.

   For $k\ge 1$,  homomorphism (\ref{wv-hom}) is the zero map.
   In fact, if $v= v^{i_1\dots i_k}
\frac{\partial}{\partial x^{i_1}}
\wedge\ldots\wedge \frac{\partial}{\partial x^{i_k}}\in {\cal V}^k(M)$,
then $v^V= v^{i_1\dots i_k}
\frac{\partial}{\partial y^{i_1}}
\wedge\ldots\wedge \frac{\partial}{\partial y^{i_k}}$.
   Denote $u^{i_1\dots i_{k-1}}=y^j \omega_{js} v^{si_1\dots i_{k-1}}$, where
$\omega=\omega_{ij} dx^i\wedge dx^j$ is the symplectic form corresponding to $w$.
   Then $u=u^{i_1\dots i_{k-1}} \frac{\partial}{\partial y^{i_1}}
\wedge\ldots\wedge \frac{\partial}{\partial y^{i_{k-1}}}$ is  a
 well-defined multivector field on~$TM$.
  One can easily verify that
$v^V=[w^V, u]=\sigma_{w^V}u$.
\end{example}

  The following example shows that the homomorphism (\ref{wv-hom})
may have a non-zero kernel in each dimension except for zero dimension.

\begin{example}
   Let, as in Example~\ref{ex-torus-Rk} $M=\T^2\times \R^k$ and
$w=\frac{\partial}{\partial x^1}\wedge \frac{\partial}{\partial x^2}$.
   Then $w^V=\frac{\partial}{\partial y^1}\wedge
\frac{\partial}{\partial y^2}$.

%

   In the dimension 1,
$$
\begin{array}{l}
H^1_P(M,w)\cong
{\cal V}^1(\R^k) \oplus C^\infty(\R^k)\oplus  C^\infty(\R^k),\\[5pt]
H^1_P(TM,w^V)\cong
{\cal V}^1(\T^2\times \R^{2k}).
\end{array}
$$
   The cohomology classes of the vertical lifts
$(f^i \frac{\partial}{\partial x^i})^V=
f^i \frac{\partial}{\partial y^i}$ vanish  and the cohomology classes of the
vertical lifts
$(f^a\frac{\partial}{\partial t^a})^V =
f^a\frac{\partial}{\partial s^a}$ are linearly independent in $H^1_P(TM,w^V)$.
   Therefore, in the dimension 1, the kernel of~(\ref{wv-hom})
is isomorphic to $C^\infty(\R^k)\oplus C^\infty(\R^k)$, and the image
of~(\ref{wv-hom}) is isomorphic to ${\cal V}^1(\R^k)$.

   In the dimensions  $\ell=2,\ldots, k+2$,
$$
\begin{array}{l}
H^\ell_P(M,w)\cong
{\cal V}^\ell(\R^k)\oplus {\cal V}^{\ell-1}(\R^k)\oplus
{\cal V}^{\ell-1}(\R^k)\oplus {\cal V}^{\ell-2}(\R^k),\\[5pt]
H^\ell_P(TM,w^V)\cong
{\cal V}^\ell(\T^2\times \R^{2k}).
\end{array}
$$
   The cohomology classes
$[(f^{a_1\dots a_\ell}\frac{\partial}{\partial t^{a_1}}\wedge \ldots
\wedge \frac{\partial}{\partial t^{a_\ell}})^V]$ are linearly independent and
generate the image  of~(\ref{wv-hom}).
   Hence, the kernel of~(\ref{wv-hom})
is isomorphic to
${\cal V}^{\ell-1}(\R^k)\oplus
{\cal V}^{\ell-1}(\R^k)\oplus {\cal V}^{\ell-2}(\R^k)$, and the image
of~(\ref{wv-hom}) is isomorphic to ${\cal V}^\ell(\R^k)$.
\end{example}

\begin{remark}
   As it seems,  there are no
natural nonzero homomorphisms between the cohomology spaces
$H^*_P(T^\A M, w^C)$ and $H^*_P(T^\A M, w^V)$.
   Simple examples show that the identity map
${\rm id}_{{\cal V}^*(T^\A M)}$ in general is not a cochain map of the complexes
$$
({\cal V}^*(T^\A M), \sigma_{w^C}) \quad \hbox{and} \quad
({\cal V}^*(T^\A M), \sigma_{w^V}).
$$
\end{remark}

\begin{proposition}
  Let $(M,w)$ be a Poisson manifold.
  For the complete lift $w^C$, the vertical lift $w^V$ of $w$ to the Weil bundle $T^\A M$,
and for any exterior form $\xi\in \Omega^*(M)$ we have

\medskip
{\rm 1)} $\widetilde{w^C}(\pi_\A^*\xi)=
(\widetilde w \xi)^V${\rm ,}

\medskip
{\rm 2)} $\widetilde{w^V}(\pi_\A^*\xi)=0${\rm ,}

\medskip
{\rm 3)} $\widetilde{w^V}(\xi^C)=\left\{
\begin{array}{ll}
 (\widetilde w\xi)^V, \quad &  \hbox{\rm if }\,|\xi|=1,\\
 0, & \hbox{\rm if }\,|\xi|\ge 2.
\end{array}
\right.$
\end{proposition}

\begin{Proof}
 1)  Let in terms of a local coordinate system
$\xi=\xi_{i_1\dots i_k} dx^{i_1}\wedge\ldots\wedge dx^{i_k}\in \Omega^k(M)$.
  Then $\pi_\A^*\xi=\xi_{i_1\dots i_k}
dx^{i_10}\wedge\ldots\wedge dx^{i_k0}$.
  It follows from Remark~\ref{cl-inj} that $(w^C)^{i0jb}=0$ for $b<n$
and  $(w^C)^{i0jn}=w^{ij}$.
  The rest of the proof is obvious  from~(\ref{vert-lift}).

2) The proof follows from~(\ref{vert-lift}).

3) For $|\xi|=1$ we have $(\xi^C)_{in}=p((\xi_i)^\A e_n) = \xi_i$ and
$\widetilde{w^V}(\xi^C)=w^{ij}\xi_i \frac{\partial}{\partial x^{jn}}=
(\widetilde w\xi)^V$.
   For $|\xi|=k\ge 2$ we have $(\xi^C)_{i_1 n i_2n \dots i_k n} =
p((\xi_i)^\A e_n\dots e_n) = 0$.
\end{Proof}

\begin{remark}
\label{weps}
  Let  $\{e_a\}$ be a Jordan-H\"older basis  in $\A$.
  From  Proposition~\ref{prop-Sch-br} it follows that
for every $a=0,1,\dots, n$ the $a$-lift
$$
w_a:=w^{(a)}=R(e_a w^\A)
$$
of $w$ is a Poisson tensor on $T^\A M$ and that these Poisson tensors are
pairwise compatible, that is,
$$
[w_a,w_b]=0.
$$

  In addition, for every $\varepsilon=\varepsilon^a e_a\in \A$ the
bivector field
\begin{equation}
\label{w-eps}
w_\varepsilon:=R(\varepsilon w^\A)
\end{equation}
on $T^\A M$ also is the Poisson tensor on~$T^\A M$ and
$w_\varepsilon=\varepsilon^a w_a$.
\end{remark}

\subsection{Modular classes of lifts of Poisson structures}
\label{subsec-mod}

   Let $(M,w)$ be an orientable Poisson manifold, and let
${\cal A}=\{(U_\kappa, h_\kappa)\}_{\kappa\in K}$ be the maximal  oriented
atlas on~$M$~\cite{Post}.
   The atlas $\cal A$ induces the  oriented     atlas
$\overline{\cal A}=
\{(\overline{U}_\kappa, \overline{h}_\kappa)\}_{\kappa\in K}$,
$\overline{U}_\kappa=\pi_\A^{-1}(U_\kappa)$,
$\overline{h}_\kappa=h^\A_\kappa$, on~$T^\A M$.
   It follows from~(\ref{ch-coord}) that the Jacobian
$\det\|\frac{\partial x^{i'}}{\partial x^{i}}\|$ of a coordinate change on~$M$  and the Jacobian
$\det\|\frac{\partial x^{i'a'}}{\partial x^{ia}}\|$
of the corresponding coordinate change on~$T^\A M$ satisfy the following relation
\begin{equation}
\label{Jac-rel}
\det\left\|
\frac{\partial x^{i'a'}}{\partial x^{ia}}\right\| =
\left(\det\left\|
\frac{\partial x^{i'}}{\partial x^{i}}\right\|\right)^{n+1}, \quad n+1=\dim\, \A.
\end{equation}

   Let $\mu$ be a volume form on $M$ and let
$$
\mu_{(U_\kappa, h_\kappa)}=\rho^{(U_\kappa, h_\kappa)}
dx^1\wedge \ldots \wedge dx^m, \quad m={\rm dim}\, M,
$$
be the coordinate representation of $\mu$.
   The family $\rho=\{\rho^{(U_\kappa, h_\kappa)}\}_{\kappa\in K}$
defines a smooth density on~$M$~\cite{Post}.
   We let
\begin{equation}
\label{orho}
\overline \rho^{(\overline{U}_\kappa,
\overline{h}_\kappa)} = ( \rho^{(U_\kappa, h_\kappa)} )^{\dim \A}.
\end{equation}
   From~(\ref{Jac-rel}) it follows that the family
$\overline \rho =
\{\overline \rho^{(\overline{U}_\kappa,
\overline{h}_\kappa)}\}_{\kappa\in K}$
defines a smooth density on $T^\A M$.
   Then the exterior form $\overline\mu$ with the coordinate representation
$$
 \overline\mu =
  \overline \rho^{(\overline{U}_\kappa, \overline{h}_\kappa)}
   dx^{10}\wedge \ldots \wedge dx^{m0}\wedge\ldots \wedge
    dx^{1n}\wedge \ldots \wedge dx^{mn},
$$
in every local chart $(\overline{U}_\kappa, \overline{h}_\kappa)$ is a volume
form on~$T^\A M$.

   Let $\Delta_{\mu}$ be the modular vector field of an oriented Poisson manifold
$(M, w, \mu)$.

  In this subsection, we compute the modular class of a Poisson structure
$w_\varepsilon$~(\ref{w-eps}) on $T^\A M$ defined by an arbitrary $\varepsilon\in\A$.
  Let $\{e_a\}$ be a Jordan-H\"older basis in $\A$, $\varepsilon=\varepsilon^a e_a$.
  We will consider  the two cases:
1) $\varepsilon^0\ne 0$, that is, $\varepsilon\notin \Ao$,
2) $\varepsilon^0=0$, that is, $\varepsilon\in \Ao$.

\begin{theorem}
\label{th-mod}
  Let $(M, w, \mu)$ be an oriented Poisson
mani\-fold, $(\A,q)$ the $(n+1)$-dimen\-sional Frobenius Weil algebra and
$\varepsilon\in \A$.

 { i)} If~$\varepsilon\notin \Ao$, then the modular vector field of
$(T^\A M, w_\varepsilon, \overline\mu)$ is
\begin{equation}
\label{mvf-we}
\Delta_{\overline \mu, w_\varepsilon} = \varepsilon^0  (n+1) \Delta_\mu^V.
\end{equation}
In particular, the modular vector field of $(T^\A M, w^C, \overline\mu)$ is
\begin{equation}
\label{mvf-l}
\Delta_{\overline \mu, w^C} = (n+1)  \Delta_\mu^V.
\end{equation}

ii)
  If~$\varepsilon\in \Ao$, then the modular vector field of
 $(T^\A M, w_\varepsilon, \overline\mu)$ is zero.
  In particular, the modular vector field of  $(T^\A M, w^V, \overline\mu)$ is zero.
\end{theorem}

\begin{Proof}
  It suffices: to verify relation~(\ref{mvf-l}) and
to~show that the modular vector field of
a Poisson structure $w_c=R(e_c w^\A)$ is zero when $c\ge 1$.
  By (\ref{mod-vf}),
in terms of a local chart $(U_\kappa, h_\kappa)$ on $M$, the modular vector
field $\Delta_\mu$ is of the form
$$
  \Delta_{\mu} =
   \sum_{j=1}^m \left(
    \frac{\partial w^{ij}}{\partial x^j} +
    w^{ij}  \frac{\partial\, \ln \rho^{(U_\kappa, h_\kappa)}}
        {\partial x^j} \right) \frac{\partial }{\partial x^i}.
$$
  From (\ref{orho}) it follows  that, in terms of the local chart
$(\overline{U}_\kappa, \overline{h}_\kappa)$ on $T^\A M$, we have
\begin{equation}
\label{dorho}
 \frac{\partial\, \ln  \overline\rho{}^{(\overline{U}_\kappa,
     \overline{h}_\kappa)}  }{\partial x^{jb}} =
  \left\{
   \begin{array}{ll}
    \displaystyle
   (n+1) \frac{\partial\, \ln \rho^{(U_\kappa, h_\kappa)} }{\partial x^{j}},
   \qquad & b=0,\\
   0, & b=1,2,\dots,n.
   \end{array}
  \right.
\end{equation}
\par
 1) First, we need to show that
\begin{equation}
\label{dwij}
 \frac{\partial (w^C)^{iajb}}{\partial x^{jb}} =
  \left\{
   \begin{array}{ll}
    \displaystyle
      \frac{\partial w^{ij}}{\partial x^{j}}, \qquad & a=n,\\[7pt]
      0, & a=0,1,\dots, n-1.
   \end{array}
  \right.
\end{equation}
  By (\ref{wC}), we have
$(w^C)^{iajb}=(w^{ij})^s \gamma^{ab}_s=
(w^{ij})^s  q^{ac}\gamma^{b}_{sc}$.

  Arguing as in Remark~\ref{rem-vl}, we find that
$$
\frac{\partial (w^{ij})^s}{\partial x^{jb}} = 0\quad \hbox{for~$s<b$},
\quad \hbox{and} \quad
\frac{\partial (w^{ij})^b}{\partial x^{jb}} =
\frac{\partial w^{ij}}{\partial x^{j}},
$$
where there is no summation over $j$ or $b$.
  Since  $\gamma_{cs}^b=0$ when $s>b$, it follows that
the only nonzero summand in
$$
\frac{\partial (w^C)^{iajb}}{\partial x^{jb}} =
\frac{\partial (w^{ij})^s  q^{ac}\gamma^{b}_{sc}}{\partial x^{jb}}
\quad
\hbox{(no summation over $j$ or $b$)}
$$
 corresponds to $s=b$.
  Since $p_n=p(e_n)=1$, it follows that
$\gamma^{ab}_b=q^{ad}\gamma^b_{bd}=1$  when $a=n$
and $\gamma^{ab}_b=0$ when $a\ne n$.
  Thus,
$$
\frac{\partial (w^C)^{iajb}}{\partial x^{jb}} = 0,
\quad a\ne n,
\quad
\hbox{and}
\quad
\frac{\partial (w^C)^{injb}}{\partial x^{jb}} =
\frac{\partial (w^{ij})^b}{\partial x^{jb}} =
\frac{\partial w^{ij}}{\partial x^{j}},
$$
where, as above, there is no summation over $j$ or $b$.
\par
  Now we show that
\begin{equation}
\label{wiaj0}
 w^{iaj0} =
  \left\{
   \begin{array}{ll}
    \displaystyle
      w^{ij}, \qquad & a= n,\\[2pt]
      0, & a=0,1,\dots, n-1.
   \end{array}
  \right.
\end{equation}
  In fact, $(w^C)^{iaj0}=(w^{ij})^s \gamma^{a0}_s$ where
$\gamma^{a0}_s = q^{ac} \gamma^0_{sc}$.
  The only nonzero summand in $q^{ac} \gamma^0_{sc}$ corresponds to $c=s=0$.
  In addition, $q^{a0}=0$ when $a\ne n$ and $q^{n0}=1$.
  We also have $(w^{ij})^0=w^{ij}$, which implies~(\ref{wiaj0}).
\par
  The modular vector field of the complete lift $w^C$ on $T^\A M$ is
$$
\displaystyle
  \Delta_{\overline\mu, w^C} =
   \sum_{jb} \left(
    \frac{\partial (w^C)^{iajb}}{\partial x^{jb}} +
    (w^C)^{iajb}
       \frac{\partial\, \ln  \overline\rho{}^{(\overline{U}_\kappa,
          \overline{h}_\kappa)}  }{\partial x^{jb}}
    \right)
    \frac{\partial }{\partial x^{ia}}.
$$
  Since the index  $b$ runs through $n+1$ different values, (\ref{dwij})
implies that
$$
\displaystyle
\sum_{jb}  \frac{\partial (w^C)^{iajb}}{\partial x^{jb}}
  \frac{\partial }{\partial x^{ia}} =
(n+1) \frac{\partial w^{ij}}{\partial x^{j}}
  \frac{\partial }{\partial x^{in}}.
$$

\par
  By (\ref{dorho}), all  summands with $b$ different from 0 in the sum
$$
\displaystyle
\sum_{jb}  (w^C)^{iajb}
       \frac{\partial\, \ln  \overline\rho{}^{(\overline{U}_\kappa,
          \overline{h}_\kappa)}  }{\partial x^{jb}}
    \frac{\partial }{\partial x^{ia}}
$$
are zero.
  Then, using (\ref{wiaj0}), we obtain
$$
\displaystyle
\sum_{jb}  (w^C)^{iajb}
       \frac{\partial\, \ln  \overline\rho{}^{(\overline{U}_\kappa,
          \overline{h}_\kappa)}  }{\partial x^{jb}}
    \frac{\partial }{\partial x^{ia}} =
   (n+1)\sum_{j}  w^{ij}
  \frac{\partial\, \ln \rho^{(U_\kappa, h_\kappa)}} {\partial x^j}
  \frac{\partial }{\partial x^{in}},
$$
which proves~(\ref{mvf-l}).
\par
  2) Consider now a Poisson structure $w_c=R(e_c w^\A)$, $c\ge 1$.
   By virtue of (\ref{pee}), we have $(w_c)^{iajb}=p((w^{ij})^s e_s e_c e^a e^b)=
(w^{ij})^s \gamma^f_{cs}\gamma^{ab}_g p(e_f e^g)= (w^{ij})^s
\gamma^f_{cs}\gamma^{ab}_f$.
  The modular vector field of $w_c$  is
\begin{equation}
\label{dmua}
\displaystyle
  \Delta_{\overline\mu, w_c} =
   \sum_{jb} \left(
    \frac{\partial (w_c)^{iajb}}{\partial x^{jb}} +
    (w_c)^{iajb}
       \frac{\partial\, \ln  \overline\rho{}^{(\overline{U}_\kappa,
          \overline{h}_\kappa)}  }{\partial x^{jb}}
    \right)
    \frac{\partial }{\partial x^{ia}}.
\end{equation}

Let us show that the first summand in the brackets in (\ref{dmua}) is zero for
all values of indices.
   By (\ref{gamma-gamma}), we have
\begin{equation}
\label{dwa-jb}
\frac{\partial (w_c)^{iajb}}{\partial x^{jb}}=
\gamma_{cs}^f\gamma^{ab}_f \frac{\partial (w^{ij})^s}{\partial x^{jb}}=
\gamma_{sf}^b\gamma^{af}_c \frac{\partial (w^{ij})^s}{\partial x^{jb}}.
\end{equation}
  Now, arguing as in the case of the structure $w^C$, we
conclude that the only nonzero summand in the right-hand side of (\ref{dwa-jb})
corresponds to $b=s$.
  But for the Jordan-H\"older basis $\gamma_{bf}^b\ne 0$
(no summation over $b$) only when $f=0$.
  On the other hand, $\gamma^{0a}_c$ does not vanish only when $a=n$ and $c=0$.
  Thus, the first summand in the brackets in (\ref{dmua}) vanishes identically.

  Let us consider  the second summand in the brackets in (\ref{dmua}).
Since $\overline\rho{}^{(\overline{U}_\kappa, \overline{h}_\kappa)}$ does not
depend on $(x^{jb})$ for $b>0$, it remains to consider only the case when
$b=0$.
  The coordinates  $(w_c)^{iaj0}=(w^{ij})^s \gamma^f_{cs}\gamma^{a0}_f$
are nonzero only when $f=0$, $a=n$.
  But $\gamma^0_{cs}=0$ if $c\ge 1$, which completes the proof.
\end{Proof}

\begin{corollary}
  In the hypotheses of Theorem~\ref{th-mod},
the modular class of the Poisson manifold $(T^\A M,w^C)$ is represented by the
vector field $(n+1) \Delta_\mu^V$ for any  modular vector field
$\Delta_\mu$  of the base manifold
$(M,w)$.
  The modular class of the Poisson manifold $(T^\A M,w^V)$
is zero.
\end{corollary}

\begin{theorem}
 Let $(M, w, \mu)$ be an oriented Poisson
mani\-fold and $(\A,q)$ the $(n+1)$-dimen\-sional Frobenius Weil algebra.
  The modular class of  $(T^\A M,w^C,\overline \mu)$ vanishes if and only if
the modular class of  $(M,w,\mu)$ vanishes.
\end{theorem}

\begin{Proof}
  Let  $\Delta_\mu$ be the modular vector
field of $(M,w)$.
  Suppose that $[\Delta_\mu]=0$, that is, $\Delta_\mu=X_g^w=[w,g]$ for some $g\in C^\infty(M)$.
  Then, by Theorem~\ref{th-mod}, the modular vector field of $(T^\A M, w^C)$ is
$\Delta_{\overline\mu}=(n+1)\Delta_\mu^V=(n+1)[w,g]^V=(n+1)[w^C,g^V]=
(n+1)[w^C,g\circ \pi_\A]$.
  Therefore $[\Delta_{\overline\mu}]=0$.

  Conversely, let $[\Delta_{\overline\mu}]=0$, then
$\Delta_{\overline\mu}=[w^C, f]$ for some $f\in C^\infty(T^\A M)$.
  Let $s_\A:M\to T^\A M$ denote  the zero section, and
let
$$
g=\frac{1}{n+1}(f\circ s_\A)\in C^\infty(M).
$$
  We claim that $\Delta_\mu=X_g^w$.

  By virtue of~(\ref{mvf-l}) and~(\ref{wC}), in terms of local coordinates,
$$
   \Delta_{\overline\mu}=
     (n+1) \sum_{k} \left(
    \frac{\partial w^{jk}}{\partial x^k} +
    w^{jk}  \frac{\partial\, \ln \rho}
        {\partial x^k} \right) \frac{\partial }{\partial x^{jn}}
$$
and
$$
    X_f^{w^C}=
       (w^C)^{iajb} \frac{\partial f}{\partial x^{ia}}
         \frac{\partial }{\partial x^{jb}}=
       (w^{ij})^s \gamma^{ab}_s \frac{\partial f}{\partial x^{ia}}
         \frac{\partial }{\partial x^{jb}},
$$
where $\mu=\rho \, dx^1\wedge\ldots \wedge dx^m$.
  From the condition $\Delta_{\overline\mu}=X_f^{w^C}$ we obtain
\begin{equation}
\label{mcl}
\begin{array}{l}
\displaystyle
(w^{ij})^s \gamma^{an}_s \frac{\partial f}{\partial x^{ia}}=
  (n+1) \sum_{k} \left(
    \frac{\partial w^{jk}}{\partial x^k} +
    w^{jk}  \frac{\partial\, \ln \rho} {\partial x^k} \right),\\[15pt]
\displaystyle
(w^{ij})^s \gamma^{ab}_s \frac{\partial f}{\partial x^{ia}} = 0 \quad \hbox{for} \quad b=0,\dots, n-1.
\end{array}
\end{equation}
  Note that
$$
  \frac{\partial f}{\partial x^{ia}}\circ s_\A=
    \frac{\partial (f\circ s_\A)}{\partial x^{ia}}.
$$
  The restriction of $(w^C)^{iajb}$ to the zero section is
$w^{ij} \gamma^{ab}_0=w^{ij} q^{ab}$.
  Thus, restricting (\ref{mcl}) to the zero section, we obtain
\begin{equation}
\label{mcl-r}
\begin{array}{l}
\displaystyle
w^{ij} q^{an} \frac{\partial (f\circ s_\A)}{\partial x^{ia}}=
  (n+1) \sum_{k} \left(
    \frac{\partial w^{jk}}{\partial x^k} +
    w^{jk}  \frac{\partial\, \ln \rho} {\partial x^k} \right),
    \\[15pt]
\displaystyle
w^{ij} q^{ab} \frac{\partial (f\circ s_\A)}{\partial x^{ia}} = 0 \quad \hbox{for} \quad b=0,\dots, n-1.
\end{array}
\end{equation}
  Therefore, contracting
the left-hand side of  (\ref{mcl-r}) with $p_b$  by virtue of~(\ref{q-p}), we
obtain
$$
w^{ij}  \frac{\partial (f\circ s_\A)}{\partial x^{ia}} q^{ab} p_b=
w^{ij}  \frac{\partial (f\circ s_\A)}{\partial x^{ia}} \delta^a_0=
w^{ij}  \frac{\partial (f\circ s_\A)}{\partial x^{i}}.
$$
  Since $p_n=1$, the contraction of the right-hand side  of (\ref{mcl-r}) with $p_b$ gives
$$
(n+1) \sum_{k} \left(
    \frac{\partial w^{jk}}{\partial x^k} +
    w^{jk}  \frac{\partial\, \ln \rho} {\partial x^k} \right).
$$
Conversely
$$
 \sum_{k} \left(
    \frac{\partial w^{jk}}{\partial x^k} +
    w^{jk}  \frac{\partial\, \ln \rho} {\partial x^k} \right) =
      w^{ij}  \frac{\partial (\frac{1}{n+1} (f\circ s_\A))}{\partial x^{i}}=
      w^{ij}  \frac{\partial g}{\partial x^{i}}.
$$

\end{Proof}

\begin{remark}
  In the case when a Poisson manifold $(M,w)$ is non-orientable,
all the  results of this subsection  remain valid.
  One only should  consider  smooth densities instead of  volume forms.
\end{remark}

{\bf Acknowledgement.}
  The author  wishes to express his deep gratitude to
Professor Mikhail A. Malakhaltsev for suggesting the problem and for many
useful conver\-sa\-tions on the subject of this paper.

%
%

%
%

\bigskip

\begin{flushleft}
\sc
Geometry Department\\
Branch of Mathematics\\
Chebotarev Research Institute of Mathematics and Mechanics\\
Kazan State University\\
Universitetskaya, 17, Kazan, 420008\\
Russia\\[5pt]
{\it E-mail:} {\tt vadimjr@ksu.ru, vshjr@yandex.ru}
\end{flushleft}

\end{document}